\documentclass[review,3p]{elsarticle}

\usepackage[english]{babel}
\usepackage{amsmath}
\usepackage{empheq}
\usepackage{amsfonts}
\usepackage{amssymb}
\usepackage{graphicx}
\usepackage{array}
\usepackage{tabularx}
\usepackage{color}
\usepackage{tcolorbox}
\usepackage{siunitx}

\usepackage{epsfig}
\newcommand{\um}{{\frac{1}{2}}}
\newcommand{\ip}{{i+\frac{1}{2}}}
\newcommand{\im}{{i-\frac{1}{2}}}
\newcommand{\QevL}{{\widetilde{\mathbf{Q}}_i^L}}
\newcommand{\QevR}{{\widetilde{\mathbf{Q}}_i^R}}

\usepackage{lineno}

\journal{Advances in Water Resources} 

\begin{document}

\begin{frontmatter}



\title{A splitting scheme for the coupled Saint-Venant-Exner model}


\author[UNITN]{A. Siviglia}
\ead{annunziato.siviglia@unitn.it}
\address[UNITN]{Laboratory of Applied Mathematics, DICAM, University of Trento, Italy}

\author[ETH]{D. Vanzo \corref{cor1}}
\ead{vanzo@vaw.baug.ethz.ch}
\address[ETH]{Laboratory of Hydraulics, Hydrology and Glaciology VAW, ETH Zurich, Switzerland
}

\author[UNITN]{E. F. Toro}
\ead{eleuterio.toro@unitn.it}

\cortext[cor1]{Corresponding author}

\begin{abstract}
We present a splitting method for the one-dimensional Saint-Venant-Exner equations used for describing the bed evolution in shallow water systems. We adapt the flux vector splitting approach  of \citet{Toro:2012} and identify one subsystem of conservative equations
(advection system) and one of non-conservative equations (pressure system), both having a very simple eigenstructure compared to the full system.
The final numerical scheme is constructed  using  a Godunov-type path-conservative scheme for the pressure system and a simple conservative Godunov method for the advection system and solved following a coupled solution strategy.
The resulting first-order accurate method is extended to second order of accuracy in space and time via the ADER approach
together with an AENO reconstruction technique. Accuracy, robustness and well-balanced properties of the resulting scheme are assessed through a carefully selected suite of testcases. The scheme is exceedingly simple, accurate and robust as the sophisticated Godunov methods. 
A distinctive feature of the novel scheme is its flexibility in the choice  of the sediment transport closure formula, which makes it particularly attractive for scientific and engineering applications.

\end{abstract}

\end{frontmatter}

   
\section{Introduction}
Nowadays numerical morphodynamic models  are used for different purposes, from answering questions about basic morphodynamic research to tackle complex engineering problems \cite{shimizu:2020,siviglia:2016}. A wide variety of river \cite{siviglia:2008,duro:2016,le:2018} and near-shore engineering problems \cite{kelly:2010,postacchini:2012} are modelled using the shallow-water approach. In this context morphodynamic investigations are often conducted using the shallow-water
equations for hydrodynamics (Saint-Venant equations \cite{saint:1871}) coupled to the equation for the bed-evolution (Exner equation \cite{exner:1925}). Both components define a coupled system of partial differential equations (PDEs) for which a conservative form does not exist,  i.e. the Saint-Venant-Exner (SVE) model.

The numerical solution of the SVE model can be obtained following two different strategies, namely decoupled and coupled. In the decoupled approach the solution is obtained solving firstly the hydrodynamic equations and assuming a fixed bed configuration and later updating the bed by using the Exner equation with the updated hydrodynamic quantities \citep{Cunge1973, Defina2003, Krishnappan1985,  Wu2004}. A clear advantage of this approach is that the governing equations 
are hyperbolic and can be expressed in conservative form. 
Despite its simplicity the decoupled approach has some numerical shortcomings.  \citet{kelly:2010} demonstrated that using a decoupled approach when simulating bore-driven sediment transport may lead to a large overestimation of the net off-shore transport in the swash zone and \citet{postacchini:2012} show that the erosion of the bed can be significantly larger  when a
swash forced by a dam-break is considered. \citet{Cordier2011} conducted numerical experiments demonstrating that the decoupled approach may fail, producing unphysical
instabilities even using a robust and well-balanced numerical scheme for shallow-water system. From a physical point of view the decoupled approach is justified when  the bed weakly interacts with the hydrodynamic waves, a condition that holds only for situations far from critical conditions \cite{Carraro:2018,devries:1965,Lyn:2002}, i.e. when the Froude number ($Fr) \ll 1$ or $\gg 1$. On the contrary, the coupled approach can be applied in all conditions at the price of having a governing system of PDEs written in non-conservative  form \cite[e.g.][]{Hudson2005a}, which requires a special numerical treatment.

While in the context of systems in conservative form the Riemann invariants and the Rankine-Hugoniot conditions provide all the necessary information to derive exact or approximate solutions of the associated Riemann problem, when dealing with system containing non-conservative products the Rankine-Hugoniot conditions across shock waves do not exist. In this case the  problem is often solved using the theory developed by \citet{Dal_Maso1995} that allows
one to set the jump conditions (and thus the concept of weak solution) in terms of a given family of paths. The degree of freedom in the choice of the family of paths is eliminated following the approach proposed by \citet{Pares:2006} who    
introduced the family of path-conservative (or path-consistent) methods to properly handle non-conservative PDEs. These approaches  have experienced a great increase of popularity in the last decade and many 
path-conservative finite-volume methods have been proposed for solving the SVE model following either a centred approach (thus not using the eigenstructure of the problem) \citep{Caleffi2007,  Canestrelli:2009,Hudson2005a} or an upwind approach (which requires a detailed knowledge of the eigenstructure) \cite[e.g.][]{ carraro:2018DOT,Castro:2008}.

In this paper we present a splitting scheme for the non-conservative SVE system of equations which is solved
following a coupled solution strategy. Our starting point is the flux vector splitting approach  of \citet{Toro:2012} (TV), first put forward for the conservative one-dimensional Euler equations. Recently the TV splitting has been successfully applied to the three-dimensional Euler equations with general equation of state \cite{toro:2015}, to the equations of magnetohydrodynamics \cite{balsara:2016} and to the Baer–Nunziato equations of compressible two-phase flow \cite{tokareva:2017}.
Our splitting identifies two separate subsystem of PDEs, the advection and the pressure systems. It differs from the original TV splitting  
in two respects, namely (i) the advection contained in the continuity equation of the Saint-Venant equations is in the pressure system and (ii) the pressure system is non-conservative. The structure of the Riemann problem associated with the pressure system being always subcritical and simple to approximate using Riemann invariants. This provides all the items required for the  evaluation of the fluctuations to be used in the update formula using the Godunov-type path-conservative method theoretically introduced in \cite{Munoz:2007} and here used for the first time in numerical applications. The advection system is hyperbolic and the numerical fluxes are obtained using a simple advection method. An attractive feature of the present method is that the sediment fluxes are contained in the advection system and are evaluated as they are described by the sediment transport formula. This means that there is no need of any differentiation as is required when the entire coefficient matrix of the SVE equations is employed for the numerical simulations \cite{carraro:2018DOT, Canestrelli:2009, Castro:2008}. 
Extension to second order is obtained through application of the ADER methodology, first introduced in  \citet{Toro:2001}
and further developed in \cite{titarevtoro,titarev:2002,dumbser:2008}. ADER has also been applied to problems governed by the non-conservative SVE equations \cite{Canestrelli:2009, Canestrelli:2010,Siviglia:2013}. Polynomial reconstruction is performed employing the AENO reconstruction procedure, an averaged variant of the popular ENO method \cite{harten:1987}, recently proposed by \citet{Toro:2020}.

The paper is structured as follows: \S \ref{sec:MathModel} briefly reviews the governing equations and the closure relationship for the sediment transport and
in \S \ref{sec:Splitting} we apply the flux splitting framework to the SVE equations and
introduce the corresponding advection and pressure systems. In \S \ref{sec:Numer_sol} we present our splitting numerical method in first order mode and in \S \ref{sec:2ndorder} we extend
it to second order of accuracy.
In \S \ref{sec:num_res}  we present numerical
results for a range of carefully selected test problems to assess both
the robustness and accuracy of the schemes proposed in this paper.
Conclusions are drawn in \S \ref{sec:conclusions}.

\section{The Saint-Venant-Exner model}
\label{sec:MathModel}
We consider the one-dimensional morphodynamic Saint-Venant-Exner model which describes the flow evolution over an erodible bed. The bed is composed of uniform sediments which are transported by the flow as bedload. 
In this section we recall the governing equations, introduce the closure relationship adopted  for its solution and then we write the system in quasi-linear form. 


\subsection{The governing equations and bedload closure relationship}
%

\begin{figure}[tbp]
\centering
\includegraphics[width=0.4\columnwidth]{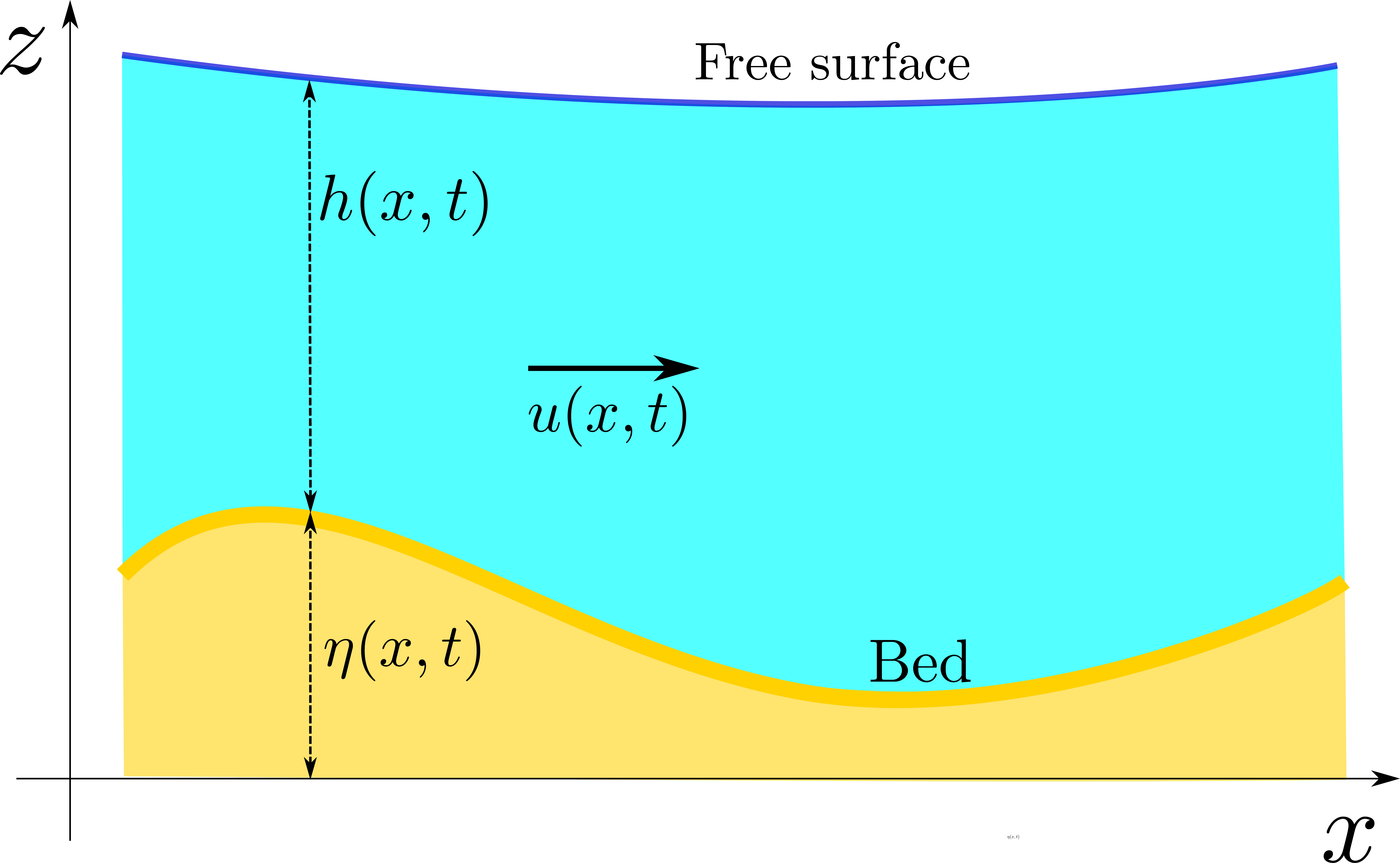}
\caption{Schematic of the water surface and erodible bottom showing notation. 
}
\label{fig:Sketch}
\end{figure}

The governing equations are obtained under shallow water conditions and includes equations for the conservation of water mass (continuity equation)
\begin{equation}\label{continuity}
\partial_t h +\partial_x q = 0\;
 \end{equation}
and  momentum of the water phase
\begin{equation}\label{momentum}
 \partial_t q + \partial_x \left(\frac{q^2}{h}+\frac{1}{2} gh^2\right) + gh \partial_x \eta=-gh S_f\; .
\end{equation}
The bed evolution is described by the sediment continuity (or Exner) equation 
\begin{equation}\label{exn}
\partial_t \eta + \partial_x q_{b} = 0\;,
\end{equation}
where, $t$[m] is time, $x$[m] is the streamwise coordinate and $g$=9.806[\si{ms^{-2}}] is the acceleration due to gravity. The quantities involved are illustrated in Fig.~\ref{fig:Sketch}. Here 
 $h$[m] is the flow depth,  $\eta$[m] is the bed level, and $u$[\si{ms^{-1}}] is depth-averaged flow velocity. The flow discharge per unit width is defined as $q=uh$[\si{m^2s^{-1}}]. 
 $q_b$[\si{m^2s^{-1}}] is the bedload sediment flux per unit width divided by $(1-\lambda_p)$ where $\lambda_p$ is the bed porosity and $S_f$[-] is the friction slope, both to be specified by an appropriate closure relationship.
  
For the sake of simplicity, we model bedload sediment flux per unit width $q_b$ in a simple form in which $q_b=q_b(u)$ only, such that \cite{Grass1981}
\begin{equation}  \label{eq:Grass}
	q_b = A_g \, u ^m \,,
\end{equation}
where $A_g$ and $m>1$ are two constant parameters.
We remark that the proposed framework can  work with any closure relationship for the bedload flux.

\subsection{Quasi-linear form of the governing equations}
The  SVE model (\ref{continuity}-\ref{exn}) is a non-conservative system \cite{Canestrelli:2009,Castro:2008} and can be written as follows 
\begin{equation}\label{completesystem_cons_non_cons}
\partial_t \mathbf{Q} + \partial_x \mathbf{F}\left(\mathbf{Q}\right) +\mathbf{B}(\mathbf{Q})\partial_x \mathbf{Q} =\mathbf{S}(\mathbf{Q})\;,
\end{equation}
where
\begin{equation} 
\mathbf{Q} = 
\left[ \begin{array}{c} 
 h \\
 q \\
 \eta 
\end{array} \right];
\qquad 
\mathbf{F} = 
\left[ \begin{array}{c} 
 q \\
 \frac{1}{2}gh^2 + q^2/h \\
 q_b 
\end{array} \right] ;
\qquad
\mathbf{B}=\left[
\begin{array}{cccc}
0& 0& 0 \\
0& 0& gh \\
0& 0& 0 
\end{array}\right];
\qquad 
\mathbf{S} = 
\left[ \begin{array}{c} 
 0 \\
 -gh S_f \\
 0 
\end{array} \right]\;.
\end{equation} 
The SVE model can also be written in quasi-linear form as
\begin{equation}\label{completesystem_non_cons}
\frac{\partial \mathbf{Q}}{\partial t}  + \mathbf{A}\left(\mathbf{Q}\right) \frac{\partial \mathbf{Q}}{\partial x} =\mathbf{S}\;,
\end{equation}
where $\mathbf{A}$ is the coefficient matrix given as
\begin{equation} \label{eq:coeff_matrix}
\mathbf{A}(\mathbf{Q})=\left[
\begin{array}{cccc}
0& 1& 0 \\
\left(c^2-u^2\right) & 2u & c^2 \\
-u\psi& \psi & 0 
\end{array}\right]\;.
\end{equation}
In the coefficient matrix $c=\sqrt{gh}$ is the celerity  and 
\begin{equation}\label{defpsi}
\psi=\frac{\partial q_b}{\partial q}
\end{equation}
is a measure of the intensity of total bedload in the flow  usually in the range $0 <\psi < \xi $ of order $\mathcal{O}\left(-2\right)$ \cite{devries:1965,Lyn:2002}. $\psi$ is obtained from differentiating the sediment transport formula and thus depends on the bedload closure relationship adopted. 
In this work, from equation ~\eqref{eq:Grass} we obtain  
\begin{equation}
\psi = m \frac{q_b}{q}.
\label{eq:phi_grass}
\end{equation}

The characteristic polynomial  of the coefficient matrix $\mathbf{A}$ is obtained by setting $|\mathbf{A} -\lambda \mathbf{I}|=0$, where $\mathbf{I}$ is the 3$\times$3-identity matrix:
\begin{equation} \label{eq_char_pol}
\lambda^3-2  u \lambda^2 
+ \left(\frac{u^2}{gh} - \psi-1\right) gh  \lambda +  ugh \psi = 0\;.
\end{equation}
If a power law formula for the solid transport is used, as that adopted in ~\eqref{eq:Grass}, the three eigenvalues 
$\lambda_1, \lambda_2, \lambda_3$ are always real, thus the governing system is always hyperbolic \cite{Cordier2011}. It is worth remarking that, under different flow conditions, either subcritical ($Fr <1$, $Fr = u / \sqrt{gh}$) or supercritical ($Fr >1$), there are always two positive and one negative eigenvalues. From a physical point of view, under sub- or supercritical conditions ($Fr < 1$ or $Fr > 1$)
the bed interacts only weakly with the water surface and  small bottom perturbations propagate at a slower pace compared with the hydrodynamic waves, whereas under near-critical conditions ($Fr \simeq 1$) the flow is close to critical conditions and the interactions between the bed and hydrodynamic waves are quite strong. 
For background on the hyperbolicity of the SVE, the eigenvalues behaviour and the physical behaviour of small bed perturbations see \cite{Cordier2011,Lyn:2002}, for example.

In this paper, we are primarily interested in the principal part of ~\eqref{completesystem_cons_non_cons} and therefore we restrict ourselves to
the homogeneous case $\mathbf{S}$($\mathbf{Q}$) = 0.

\section{Splitting framework} 
\label{sec:Splitting}
In this section, we propose a splitting method for systems written in non-conservative form 
following the  framework of \citet{Toro:2012}. We remark that this novel formulation is valid either 
under sub- or supercritical conditions ($Fr < 1$ or $Fr > 1$) or under near-critical conditions ($Fr \simeq 1$). 
%
%

\subsection{The framework}
Consider now the homogeneous SVE equations
\begin{equation}\label{completesystem_cons_non_cons_no_source}
\partial_t \mathbf{Q} + \partial_x \mathbf{F}\left(\mathbf{Q}\right) +\mathbf{B}(\mathbf{Q})\partial_x \mathbf{Q} =\mathbf{0}
\;,
\end{equation}
with
\begin{equation} 
\mathbf{Q} = 
\left[ \begin{array}{c} 
 h \\
 q \\
 \eta 
\end{array} \right];
\qquad 
\mathbf{F}(\mathbf{Q}) = 
\left[ \begin{array}{c} 
 q \\
 \frac{1}{2}gh^2 + q^2/h \\
 q_b 
\end{array} \right] ;
\qquad
\mathbf{B}(\mathbf{Q})=\left[
\begin{array}{cccc}
0& 0& 0 \\
0& 0& gh \\
0& 0& 0 
\end{array}\right] \;.
\end{equation} 
First, we identify the conservative part and express the conservative flux as the sum of
advection and pressure fluxes as follows
\begin{equation} 
\mathbf{F}(\mathbf{Q}) = 
\left[ \begin{array}{c} 
 0 \\
 q^2/h \\
 q_b 
\end{array} \right] + 
\left[ \begin{array}{c} 
 q \\
 \frac{1}{2}gh^2 \\
 0 
\end{array} \right] 
\end{equation} 
with the corresponding advection and pressure fluxes defined as
\begin{equation} \label{eq:flux_splitting}
\mathbf{F}^{(a)}(\mathbf{Q}) = q 
\left[ \begin{array}{c} 
 0 \\
 q/h \\
 q_b/q  
\end{array} \right] \qquad \text{and} \qquad
 \mathbf{F}^{(p)}(\mathbf{Q}) =
\left[ \begin{array}{c} 
 q \\
 \frac{1}{2}gh^2 \\
 0 
\end{array} \right] \;.
\end{equation} 

Then we consider two subsystems:
\begin{empheq}[left=\empheqlbrace]{align}
  &\partial_t \mathbf{Q} + \partial_x \mathbf{F}^{(a)}(\mathbf{Q}) =\mathbf{0} \label{eq:adv} \\
  &\partial_t \mathbf{Q} +  \partial_x \mathbf{F}^{(p)}(\mathbf{Q}) + \mathbf{B}(\mathbf{Q}) \partial_x \mathbf{Q} =\mathbf{0} \label{eq:prex}
\end{empheq}
called respectively  the \textit{advection system} (\ref{eq:adv}) and the \textit{pressure system} (\ref{eq:prex}). We note however that here the pressure system is augmented by the non-conservative term present in the SVE equations.
The final goal of this procedure is to obtain the numerical solution of the full SVE system of equations. The TV flux splitting
approach consists of approximating the numerical fluxes for the pressure system and advection system separately and
constructing the numerical fluctuations for the full system based on these. To this end, the analysis of the eigenstructure and the
study of the Riemann problem for the pressure system are necessary.

\subsubsection{The pressure system} 
The pressure system (\ref{eq:prex}) is non-conservative because of the presence of the non-conservative term $g h \partial_x \eta$ appearing in the momentum equation (\ref{momentum}). It can be written in quasi-linear form  as
\begin{equation}
 \partial_t \mathbf{Q} +  \mathbf{P}(\mathbf{Q}) \partial_x \mathbf{Q} =\mathbf{0}
\end{equation}
with 
\begin{equation} \label{eq:pressure_matrix}
\mathbf{P}= \mathbf{J}^{(p)} + \mathbf{B}   
 =\left[
\begin{array}{cccc}
0& 1& 0 \\
c^2 & 0 & 0 \\
0 & 0 & 0 
\end{array}\right] +
\left[
\begin{array}{cccc}
0& 0& 0 \\
0 & 0 & c^2 \\
0 & 0 & 0 
\end{array}\right] =
\left[
\begin{array}{cccc}
0& 1& 0 \\
c^2 & 0 & c^2 \\
0 & 0 & 0 
\end{array}\right]
\end{equation} 
where $\mathbf{J}^{(p)}$ is the Jacobian matrix of the pressure fluxes $\mathbf{F}^{(p)}$ in (\ref{eq:flux_splitting}).
The eigenvalues of matrix $\mathbf{P}$ are 
\begin{equation} \label{eq:eigenvalues_P}
\lambda_{1}^{(p)} = -c,  \qquad \lambda_{2}^{(p)} = 0,  \qquad \lambda_{3}^{(p)} = c \;.
\end{equation}
The eigenvalues are always real and  $\lambda_{1}^{(p)} < \lambda_{2}^{(p)} = 0 < \lambda_{3}^{(p)}$ and thus the system is always subcritical as illustrated in Fig.~\ref{fig:RP_prex}. 
The right eigenvectors corresponding to the three eigenvalues (\ref{eq:eigenvalues_P}) are 
\begin{equation}
\label{eq:eigenvectors_P}
 \mathbf{R}_{1}^{(p)} = 
\left[ \begin{array}{r} 
 1 \\
 -c \\
 0 
\end{array} \right], 
\qquad 
\mathbf{R}_{2}^{(p)} =
\left[ \begin{array}{r} 
 -1 \\
 0 \\
 1 
\end{array} \right], 
\qquad
\mathbf{R}_{3}^{(p)} =
\left[ \begin{array}{r} 
 1 \\
 c \\
 0 
\end{array} \right]. 
\end{equation}

\subsubsection{The advection system} 
The advection system in conservative form is
\begin{equation}
\partial_t \mathbf{Q} + \partial_x \mathbf{F}^{(a)}(\mathbf{Q}) =\mathbf{0},
\end{equation}
where $\mathbf{Q}=[h,q,\eta]^T$ and $\mathbf{F}^{(a)}$ as in (\ref{eq:flux_splitting}). The quasi-linear form is given by 
\begin{equation}
\partial_t \mathbf{Q} +  \mathbf{J}^{(a)}(\mathbf{Q}) \partial_x \mathbf{Q} =\mathbf{0}, 
\end{equation} 
where 
\begin{equation} \label{eq:adv_matrix_fluxes}
\mathbf{J}^{(a)}=
\left[
\begin{array}{cccc}
0& 0& 0 \\
-u^2 & 2u & 0 \\
-u\psi& \psi & 0 
\end{array}
\right]\; 
\end{equation}
is the Jacobian of matrix $\mathbf{F}^{(a)}$. Simple analysis shows that the eigenvalues of matrix (\ref{eq:adv_matrix_fluxes}) are $\lambda_1^{(a)}=0$ and $\lambda_2^{(a)}=\lambda_3^{(a)}=u$ and that there are only two  linearly independent right eigenvectors given by
\begin{equation}
\label{eq:eigenvectors_A}
 \mathbf{R}_{1}^{(a)} = \alpha_1
\left[ \begin{array}{r} 
 0 \\
 0 \\
 1 
\end{array} \right], 
\qquad 
\mathbf{R}_{2}^{(a)} = \alpha_2
\left[ \begin{array}{r} 
 1 \\
 u \\
 0 
\end{array} \right],  
\end{equation}
thus the advection system is weakly hyperbolic. It is easy to show that the $\lambda_1^{(a)}$-field is linearly degenerate while the $\lambda_2^{(a)}$ and $\lambda_3^{(a)}$ are genuinely  non-linear if $\alpha_2 \neq 0$ and $u \neq 0$.
We note that the weakly hyperbolic nature of the advection system does not have a bearing on its numerical approximation.

\section{Numerical solution}
\label{sec:Numer_sol}
Direct integration of (\ref{completesystem_cons_non_cons_no_source}) in the space-time control volume $V_i=[x_{i-\frac{1}{2}};x_{i+\frac{1}{2}}] \times [t^n;t^{n+1}]$ gives the following 
update numerical formula to solve (\ref{completesystem_cons_non_cons_no_source}):
\begin{equation} \label{eq:update}
\mathbf{Q}_i^{n+1} = \mathbf{Q}_i^{n} - \frac{\Delta t}{\Delta x}\left[ \left( \mathbf{D}_{\ip}^- + \mathbf{D}_{\im}^+ \right) + \left(\mathbf{F}^{(a)}_{\ip} -\mathbf{F}^{(a)}_{\im} \right)   \right].   
\end{equation} 
$\mathbf{Q}_i^{n+1}$ and $\mathbf{Q}_i^{n}$ are the cell-averaged values of the vector $\mathbf{Q}$ at
time $t^{n+1}$ and $t^{n}$, $\Delta t$ is the time step derived from the standard CFL stability condition and $\Delta x$ is the grid size (here, for simplicity, assumed having a constant size).  
$\mathbf{D}_{\ip}^-$ and $\mathbf{D}_{\im}^+$ are \textit{fluctuations}, or \textit{increments} associated to the pressure system, which are obtained using path-conservative schemes \cite{Pares:2006}, while $\mathbf{F}^{(a)}_{\ip}$ and $\mathbf{F}^{(a)}_{\im}$ are the numerical fluxes of the advection system obtained using simple advection methods. 

In order to compute the fluctuations $\mathbf{D}_{\ip}^-$ and $\mathbf{D}_{\im}^+$ and  the advection fluxes $\mathbf{F}^{(a)}_{i + \um}$ and $\mathbf{F}^{(a)}_{i - \um}$  to be used in  ($\ref{eq:update}$) we consider the Riemann problem for each system. We start
with the pressure system.

\subsection{The pressure system}
\label{sec:pressure_system}
\begin{figure}[tbp]
\centering
\includegraphics[width=0.4\columnwidth]{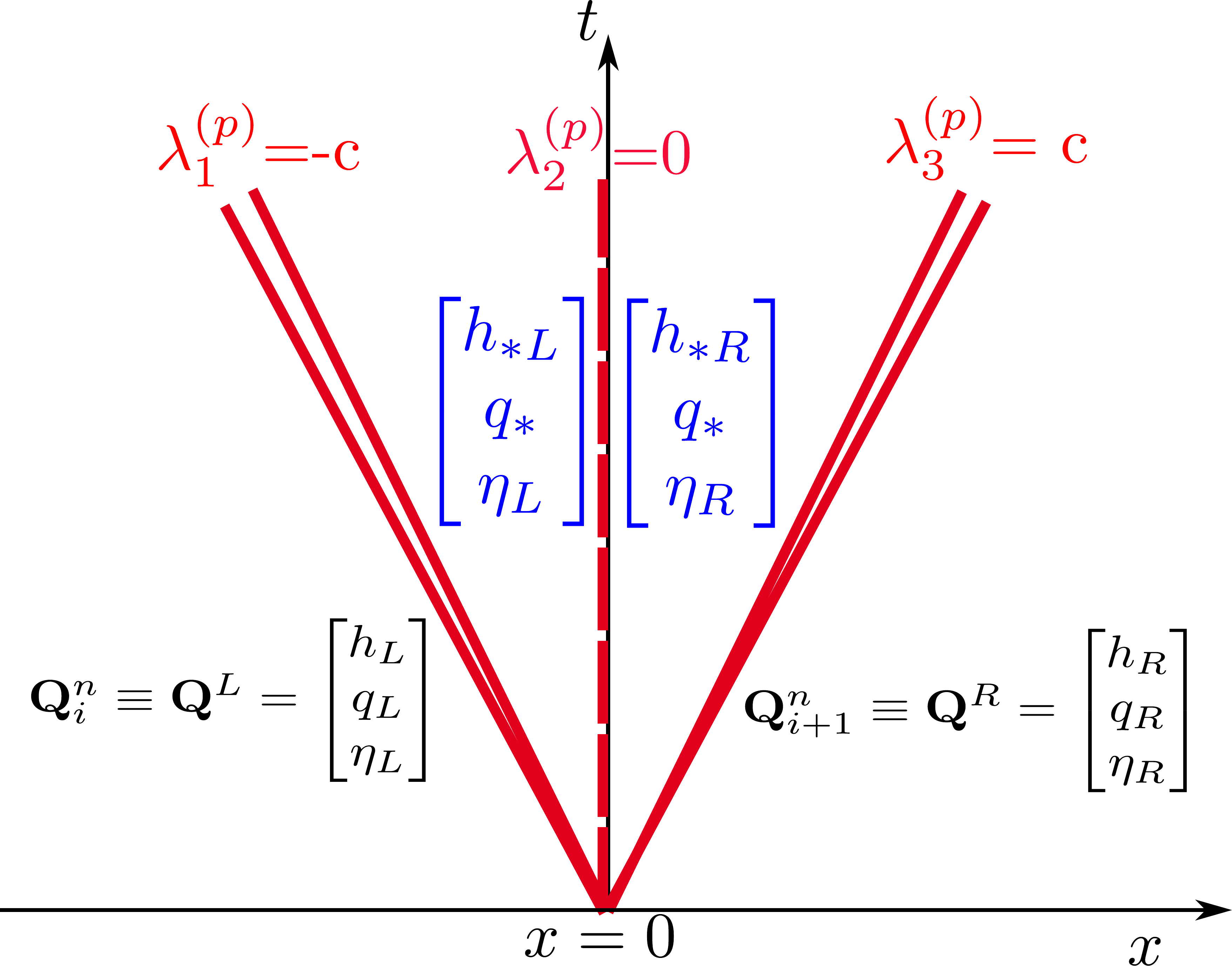}
\caption{ \textbf{Structure of the solution in local coordinates of the Riemann problem for the pressure system resulting from the flux splitting}. There are two non-linear wave
families and a stationary contact discontinuity coinciding
with the $t$-axis. The wave pattern is always subcritical therefore determining the Godunov state for flux evaluation does not require sampling, being always the star state. The sought values in the star region are  $h_{*L}$, $h_{*R}$,   $q_*$ and $\eta_{*L}$, $\eta_{*R}$. 
}
\label{fig:RP_prex}
\end{figure}

In order to calculate  the fluctuations $\mathbf{D}_{\ip}^-$ and $\mathbf{D}_{\ip}^+$ at the interface $x_{\ip}$we consider the Riemann problem for the pressure system in conservative variables 
\begin{equation}                                                 \label{eq:RP_prex}
       \left.
       \begin{array}{l}
       \partial_{t}\mathbf{Q} + \mathbf{P}(\mathbf{Q}) \partial_{x}\mathbf{Q} = 
       {\bf 0} \;,\hspace{2mm} x \in {\cal{R}}  \;, \hspace{2mm} t>0 \;, \\
       \\
       \mathbf{Q}(x,0) = \left\{ 
                      \begin{array}{lll}
                            \mathbf{Q}^{L} \equiv \mathbf{Q}_{i}^n & \mbox{ if } & x < 0  \;,\\
                            \\
                            \mathbf{Q}^{R} \equiv \mathbf{Q}_{i+1}^n & \mbox{ if } & x > 0  \;.
                      \end{array}
                      \right.
      \end{array}
      \right\}
\end{equation}
The structure of the solution of (\ref{eq:RP_prex}) at the fixed interface position $x_{\ip}$ , or $x=0$ in local coordinates, is illustrated in Fig.~\ref{fig:RP_prex}.
The wave pattern is always subcritical and composed by three wave families. The left family is associated
with the eigenvalue $\lambda_1^{(p)}$, the middle family is superimposed onto the $t$-axis, and is
associated with $\lambda_2^{(p)}$ and the right family is associated  with $\lambda_3^{(p)}$.
The waves associated with the genuinely non-linear characteristic fields $\lambda_1^{(p)}$ and $\lambda_3^{(p)}$
are either shocks (discontinuous solutions) or rarefactions (smooth solutions), while the
wave associated with the linearly degenerate characteristic field $\lambda_2^{(p)}$ is a stationary contact
discontinuity.
The entire solution consists of four constant states, namely $\mathbf{Q}^L=\left[h_L,q_L,\eta_L\right]^T$ (data), $\mathbf{Q}_{*}^L=\left[h_{*L},q_{*L},\eta_{*L}\right]^T$, $\mathbf{Q}_{*}^R=\left[h_{*R},q_{*R},\eta_{*R}\right]^T$ and  $\mathbf{Q}^{R}=\left[h_{R},q_{R},\eta_{R}\right]^T$ (data), separated by the three distinct waves. The unknown states to be found in the star region are 
$\mathbf{Q}_{*}^L$ (left of $x = 0$) and $\mathbf{Q}_{*}^R$ (right of $x = 0$). 
We apply Riemann invariants to find these solutions. 
Across the stationary contact discontinuity with right eigenvectors (\ref{eq:eigenvectors_P}) the generalized Riemann
invariants are solutions of the following two ordinary differential equations (ODEs) 
\begin{equation} \label{eq:contact}
\textit{Stationary contact} \Rightarrow \frac{d h}{-1} = \frac{d q}{0} = \frac{d \eta}{1} \; 
\end{equation}
which can be obviously rewritten as a system of ODEs
\begin{empheq}[left=\empheqlbrace]{align}
  &\frac{d h}{-1} = \frac{d q}{0} \label{eq:RI_contact_1} \\
  &\frac{d h}{-1} =  \frac{d \eta}{1} \label{eq:RI_contact_2} \;.
\end{empheq}
Integration of the \eqref{eq:RI_contact_1} in the phase space gives $q=const$  across the wave and thus $q_{*L}=q_{*R}=q_{*}$, while integration of \eqref{eq:RI_contact_2} gives
\begin{equation} \label{eq:H_stat_contact}
h_{*L}+\eta_{*L}=h_{*R}+\eta_{*R}=H_* \;.
\end{equation}
Across the left and right waves we have:
\begin{empheq}{align}
   \textit{Left wave} \Rightarrow & \quad \frac{d h}{1} = \frac{d q}{-c} = \frac{d \eta}{0}  \label{eq:first_eq_left} \\
  \textit{Right wave} \Rightarrow & \quad \frac{d h}{1} = \frac{d q}{c} = \frac{d \eta}{0}\label{eq:first_eq_right} \;.
\end{empheq}
The third ODEs in Eqs.~(\ref{eq:first_eq_left}) and (\ref{eq:first_eq_right}) imply that $\eta$ remains constant across the left (i.e. $\eta_{*L}=\eta_L$) and right wave (i.e. $\eta_{*R}=\eta_R$).
Exact integration of the first ODEs in Eqs.~(\ref{eq:first_eq_left}) and (\ref{eq:first_eq_right}) respectively gives that  
\begin{empheq}[left=\empheqlbrace]{align}
  & \frac{2}{3} \sqrt{g} h^{3/2} + q = const \qquad \text{across the left wave} \label{eq:Riemann_inv_left}\\ 
& \frac{2}{3} \sqrt{g} h^{3/2} - q = const \qquad \text{across the right wave} \label{eq:Riemann_inv_right} \;.
\end{empheq}

After simple algebraic manipulations of Eqs.~(\ref{eq:H_stat_contact}), (\ref{eq:Riemann_inv_left}) and (\ref{eq:Riemann_inv_right}), we obtain the following non-linear system 
\begin{empheq}[left=\empheqlbrace]{align}
& h_{*L}^{3/2} + h_{*R}^{3/2} = K \label{eq_non_lin_sys_h1}\\ 
& h_{*L} - h_{*R} = \Delta \eta \label{eq_non_lin_sys_h2} \\
& q_{*} = \frac{1}{2}\left(q_L+q_R\right) +\frac{\sqrt{g}}{3} \left(h_L^{3/2} -h_R^{3/2} - h_{*L}^{3/2} + h_{*R}^{3/2} \right) \label{eq_non_lin_sys_q}
\end{empheq}

where
\begin{equation} \label{eq:K_delta_eta}
K= \frac{3}{2\sqrt{g}} \left( q_L-q_R \right) +  h_L^{3/2} + h_R^{3/2},  \qquad \qquad
\Delta \eta = \eta_R -\eta_L \;.
\end{equation}
Iterative solution of (\ref{eq_non_lin_sys_h1}) and (\ref{eq_non_lin_sys_h2}) gives the sought values $h_{*L}$ and $h_{*R}$. Substitution of such values in (\ref{eq_non_lin_sys_q}) provides $q_*$.  

One could improve efficiency by finding an approximate solution of (\ref{eq_non_lin_sys_h1}) in closed form.
A possible way to avoid the iterations is the following. We assume that $\Delta \eta =0$ and solve (\ref{eq_non_lin_sys_h2}),  that gives $h_{*L}=h_{*R}= \widehat{h}$.
Then solving (\ref{eq_non_lin_sys_h1}) we obtain a closed form solution for $\widehat{h}$, i.e.
\begin{equation} \label{eq:hsL_closed}
\widehat{h}  =\left[ \frac{3}{4\sqrt{g}} \left( q_L-q_R \right)  +  \frac{1}{2} \left(h_L^{3/2} + h_R^{3/2}\right)\right]^{2/3} \;. 
\end{equation} 
At this stage we linearize (\ref{eq_non_lin_sys_h1}), which, together with  (\ref{eq_non_lin_sys_h2}) gives the following linear system  
\begin{equation} \label{eq:non_lin_system_linearized}
\left\{
\begin{array}{ll}
 h_{*L}\sqrt{\widehat{h}} + h_{*R}\sqrt{\widehat{h}}  = K \\ 
 h_{*L} - h_{*R} = \Delta \eta 
\end{array}
\right.
\end{equation}
that gives 
\begin{equation} \label{eq:hstar_linear}
h_{*L} = \frac{1}{2} \left(  \frac{K}{\sqrt{\widehat{h}}} + \Delta \eta\right) \\ \; .
\end{equation}
Finally we obtain the following system that can be directly solved: 
\begin{equation} \label{eq:non_lin_system_linearized_final}
\left\{
\begin{array}{ll}
 h_{*L} = \frac{1}{2} \left(  \frac{K}{\sqrt{\widehat{h}}} + \Delta \eta\right) \\ 
 h_{*R} = h_{*L} - \Delta \eta \\
q_{*} = \frac{1}{2}\left(q_L+q_R\right) +\frac{\sqrt{g}}{3} \left(h_L^{3/2} -h_R^{3/2} - h_{*L}^{3/2} + h_{*R}^{3/2} \right) \;.
\end{array}
\right.
\end{equation}

Once the solution in the star region is known, the fluctuations $\mathbf{D}_{\ip}^-$ and $\mathbf{D}_{\ip}^+$  are obtained using a  Godunov-type path-conservative method as in \cite{Munoz:2007}:
\begin{equation} \label{eq:fluctuations_God}
\left.
\begin{array}{ll}
\mathbf{D}_{\ip}^- = \frac{1}{\Delta t} \displaystyle\int_{0}^{\Delta t} \displaystyle\int_{0}^{1} \mathbf{P}(\mathbf{\Psi}(s; \mathbf{Q}_{i}^n,\mathbf{Q}_{*}^R)) \frac{ \partial}{\partial{s}}\mathbf{\Psi}(s; \mathbf{Q}_{i}^n,\mathbf{Q}_{*}^R)) ds \, dt \\
\mathbf{D}_{\ip}^+ = \frac{1}{\Delta t} \displaystyle\int_{0}^{\Delta t} \displaystyle\int_{0}^{1} \mathbf{P}(\mathbf{\Psi}(s; \mathbf{Q}_{*}^L,\mathbf{Q}_{i+1}^n)) \frac{ \partial}{\partial{s}}\mathbf{\Psi}(s; \mathbf{Q}_{*}^L,\mathbf{Q}_{i+1}^n)) ds \, dt\;.
\end{array}
\right\}
\end{equation} 
We remark that, in the Riemann problem associated with the pressure system,  we have two star regions, one on the left of $x=0$ (whose solution is 
$\mathbf{Q}_{*}^L$)  and one on the right (whose solution is $\mathbf{Q}_{*}^R$). The presence of these two regions is taken into account integrating the left fluctuation  $\mathbf{D}_{\ip}^-$ between $\mathbf{Q}_{i}^n$ and $\mathbf{Q}_{*}^R$ and the right fluctuation $\mathbf{D}_{\ip}^+$ between $\mathbf{Q}_{*}^L$ and $\mathbf{Q}_{i+1}^n$. We verified numerically that the fluctuations calculated as in \eqref{eq:fluctuations_God} satisfy the compatibility condition 
\begin{equation} \label{eq:compatibility}
\mathbf{D}_{\ip}^- + \mathbf{D}_{\ip}^+ = \frac{1}{\Delta t} \displaystyle\int_{0}^{\Delta t} \displaystyle\int_{0}^{1} \mathbf{P}(\mathbf{\Psi}(s; \mathbf{Q}_{i}^n,\mathbf{Q}_{i+1}) \frac{ \partial}{\partial{s}}\mathbf{\Psi}(s; \mathbf{Q}_{i}^n,\mathbf{Q}_{i+1}) ds .
\end{equation}

For all numerical test cases presented in this paper, we always use the simple 
segment paths, given by
\begin{equation} \label{eq:seg_path}
\left.
\begin{array}{lll}
 \mathbf{\Psi}(s;\mathbf{Q}^n_i,\mathbf{Q}_{*}^R) &= \mathbf{Q}^n_i &+ s \left( \mathbf{Q}_{*}^R - \mathbf{Q}^n_i \right) \\ 
 \mathbf{\Psi}(s;\mathbf{Q}_{*}^L,\mathbf{Q}^n_{i+1}) &= \mathbf{Q}_{*}^L &+ s \left( \mathbf{Q}^n_{i+1} - \mathbf{Q}_{*}^L \right). 
\end{array}
\right\}  
\end{equation}  
Then, from ~(\ref{eq:fluctuations_God}) we have
\begin{equation} \label{eq:fluctuations_God_num}
\mathbf{D}_{\ip}^- = \mathbf{\widehat{P}}_\ip^- \left[ \mathbf{Q}_{*}^R - \mathbf{Q}^n_i \right]  \qquad ; \qquad 
\mathbf{D}_{\ip}^+ = \mathbf{\widehat{P}}_\ip^+ \left[\mathbf{Q}^n_{i+1} -\mathbf{Q}_{*}^L  \right] 
\end{equation}
where  
\begin{equation} 
 \label{eq:fluctuations_God_approx} 
  \mathbf{\widehat{P}}_\ip^- \approx \frac{1}{\Delta t} \displaystyle\int_0^{\Delta t} \displaystyle\int_{0}^{1} \mathbf{P}(\mathbf{\Psi}(s; \mathbf{Q}_{i}^n,\mathbf{Q}_{*}^R)) ds \, dt \qquad ; \qquad 
\mathbf{\widehat{P}}_\ip^+ \approx \frac{1}{\Delta t} \displaystyle\int_0^{\Delta t} \displaystyle\int_{0}^{1} \mathbf{P}(\mathbf{\Psi}(s; \mathbf{Q}_{*}^L,\mathbf{Q}_{i+1}^n)) ds \, dt \;.
\end{equation} 

Given a $nGP$-point Gaussian quadrature rule with weights $\omega_j$ and positions $s_j$ distributed in the unit interval 
$[0;1]$, a very accurate numerical approximation of the matrices $\mathbf{\widehat{P}}_\ip^-$ and $\mathbf{\widehat{P}}_\ip^+$ is given by 
\begin{equation} 
 \label{eqn:integrat_Gauss}
  \mathbf{\widehat{P}}_\ip^- = \sum_{j=1}^{nGP} \omega_j \mathbf{P}(\mathbf{\Psi}(s_j,\mathbf{Q}_i,\mathbf{Q}_{*}^R)) \qquad ; \qquad 
  \mathbf{\widehat{P}}_\ip^+ = \sum_{j=1}^{nGP} \omega_j \mathbf{P}(\mathbf{\Psi}(s_j,\mathbf{Q}_{*}^L,\mathbf{Q}_{i+1})) \;. 
\end{equation}
As an example using a three-point Gaussian quadrature rule with the following points $s_j$ and 
weights $\omega_j$:
\begin{equation}
  s_1 = \frac{1}{2}, \quad s_{2,3} = \frac{1}{2} \pm \frac{\sqrt{15}}{10}, \qquad \omega_1 = \frac{8}{18}, \quad \omega_{2,3} =  \frac{5}{18}.
\end{equation}
All simulation of this work  are calculated using  $nGP$=1. This is enough to ensure the achievement of second order of accuracy.

At this stage we anticipate that there is no visible difference between the numerical results obtained with  numerical fluctuations (\ref{eq:fluctuations_God}) that uses the linearized solution (\ref{eq:hstar_linear}) and those obtained from a fluctuation that uses the iterative solution of the non-linear system (\ref{eq_non_lin_sys_h1}-\ref{eq_non_lin_sys_q}). 
We therefore recommend the linearized solution (\ref{eq:hstar_linear}) for practical applications.


\subsection{The advection system} 
Recall that in our splitting (\ref{eq:flux_splitting}) the advection operator is written as
\begin{equation}
\label{eq:Flux_adv}
 \mathbf{F}^{(a)} \left(\mathbf{Q}\right)=  
 \left[ \begin{array}{c} 
 0 \\
  q^2/h\\
 A_g (q/h)^m  
\end{array} \right] \;.
\end{equation}


The algorithm we propose for constructing the numerical flux $\mathbf{F}^{(a)}_{\ip}$ to be used in the update formula (\ref{eq:update}) is


\begin{equation}                         \label{eq:advection_flux}
     {\bf F}_{i+\frac{1}{2}}^{(a)} = \left\{\begin{array}{c}
     q_{*} \begin{bmatrix}
     0\\
     \displaystyle{(u_i^n)} \\
     A_g \frac{\left(u_i^n\right)^{m-1}}{(h_i^n)} 
     \end{bmatrix}  \mbox{ if } q_{*} \ge 0 \;, \\
     \\
     q_{*} \begin{bmatrix}
     0\\
     \displaystyle{(u_{i+1}^n)} \\
     A_g \frac{\left(u_{i+1}^n\right)^{m-1}}{(h_{i+1}^n)} 
     \end{bmatrix}  \mbox{ if } q_{*} < 0 
     \end{array} \right. 
\end{equation}
where $q_*$ is the solution (\ref{eq_non_lin_sys_q}) emerging from the Riemann problem of the pressure system.

\subsection{Summary of the proposed scheme}
\label{sec:summary}
In order to compute  the fluctuations $\mathbf{D}_{\ip}^-$ and $\mathbf{D}_{\im}^+$ and advection fluxes $\mathbf{F}^{(a)}_{i + \um}$ and  $\mathbf{F}^{(a)}_{i - \um}$  to be used in the update formula ($\ref{eq:update}$) we proceed as follows:
\begin{itemize}
\item \textbf{Pressure fluctuations}. At each interface evaluate the solution of the Riemann problem $\mathbf{Q}_{*}^L=[h_{*L},q_*,\eta_L]^T$ and $\mathbf{Q}_{*}^R=[h_{*R},q_*,\eta_R]^T$ using Eqs.~(\ref{eq:hstar_linear}), (\ref{eq_non_lin_sys_h2}) and (\ref{eq_non_lin_sys_q}). Then calculate  $\mathbf{D}_{\ip}^-$ and $\mathbf{D}_{\ip}^+$ as in (\ref{eq:fluctuations_God_num}) evaluating the approximated matrices (\ref{eqn:integrat_Gauss}) using one Gaussian point.
\item \textbf{Advection flux} Evaluate the advection fluxes  $\mathbf{F}^{(a)}_{\ip}$ 
as described in (\ref{eq:advection_flux}).
\end{itemize}

\section{Second order extension}
\label{sec:2ndorder}
Extension to second order is obtained using the ADER methodology  by \citet{Toro:2001}. 
The procedure to achieve second order contains two ingredients: (i) a first-order non-linear spatial reconstruction of the gradient of the solution in each cell and (2) the solution of the generalized Riemann problem (GRP) at the interface of each cell. For background on ADER see Chapters 19 and 20 of of \cite{Toro:2013} and references therein.  
Here we consider second-order accurate ADER schemes based on the HEOC
solver of \citet{harten:1987} (see also \cite{castro_crist:2008}) for the GRP.

\subsection{Nonlinear Reconstruction Technique}
\label{sec.rec.nonlinear}
First we deal with the reconstruction
problem. We adopt the AENO reconstruction procedure, an averaged variant of the popular ENO method \cite{harten:1987}, recently proposed by \citet{Toro:2020}. To achieve second-order of accuracy we need to construct first-degree
polynomials $\boldsymbol{p}_i$ in each cell $I_i$ at time 
$t^n$ from the given cell averages $\left\{\mathbf{Q}_i^n\right\}$  
of the form
\begin{equation} \label{eq:polynm_space}
\boldsymbol{p}_i = \mathbf{Q}^n_i + (x-x_i) \mathbf{\Delta}_i
\end{equation}
where $\mathbf{\Delta}_i$ is the slope vector and $x_i=(x_\im+x_\ip)/2$.
Recall that in order to circumvent Godunov’s Theorem \cite{Godunov:1959}, the reconstruction must be non-linear.
See Chapter 20 of \cite{Toro:2013} for background. Here the non-linearity of the scheme is ensured by taking the polynomial slope as
\begin{equation} \label{eq:slope_AENO}
\mathbf{\Delta}_i = \um(1+\beta) \mathbf{\Delta}_\im  + \um(1-\beta) \mathbf{\Delta}_\ip \qquad \text{with} \qquad |\beta| \le 1
\end{equation}
where
\begin{equation} \label{eq:omega_AENO}
\beta(r) = \frac{1-r}{\sqrt{\epsilon^2+(r-1)^2}} \qquad \text{with} \qquad r = \frac{\left| \mathbf{\Delta}_\im \right| }{\left| \mathbf{\Delta}_\ip \right| + TOL} 
\end{equation}
and 
\begin{equation} \label{eq:DELTA_AENO}
\mathbf{\Delta}_\im = \frac{\mathbf{Q}^n_i - \mathbf{Q}^n_{i-1} }{\Delta x} \qquad , \qquad \mathbf{\Delta}_\ip = \frac{\mathbf{Q}^n_{i+1} - \mathbf{Q}^n_{i} }{\Delta x} \;.
\end{equation}
The parameter $\epsilon$ is a positive constant, while $TOL$ is a small positive tolerance to avoid division by zero. 

\subsection{ Second-order ADER with the HEOC solver for the GRP}
\label{sec.rec.timedisc}
The result of the reconstruction procedure is a non-oscillatory linear polynomial $\boldsymbol{p}_i$  defined at time $t^n$ inside each 
spatial element $I_i$. 
We are interested in the left and right limiting values of the reconstruction polynomials, often called
boundary extrapolated values.
Let us first consider cell $I_i$ with cell boundaries $x_\im$ and $x_\ip$ and define
\begin{equation} \label{eq:boundary_values}
\left.
\begin{array}{ll}
 \mathbf{Q}_i^L = \boldsymbol{p}_i(x_\im)= \mathbf{Q}(x_\im^+,0), \\ 
\mathbf{Q}_i^R = \boldsymbol{p}_i(x_\ip)= \mathbf{Q}(x_\ip^-,0) \;. 
\end{array}
\right\}
\end{equation}
We now evolve these limiting values in time using the time Taylor series expansion
\begin{equation} \label{eq:boundary_evolved}
\left.
\begin{array}{ll}
 \QevL(\tau) = \mathbf{Q}(x_\im^+,0) + \tau \partial_{t}\mathbf{Q}(x_\im^+,0), \\ 
 \QevR(\tau) = \mathbf{Q}(x_\ip^-,0) + \tau \partial_{t}\mathbf{Q}(x_\ip^-,0) \;. 
\end{array}
\right\}
\end{equation}
After adopting notation (\ref{eq:boundary_values}) and using the Cauchy-Kovalevskaya procedure,  in the case of our splitting, the time derivatives above can be expressed as
\begin{equation} \label{eq:Kovalevskaya}
\left.
\begin{array}{ll}
 \partial_{t}\mathbf{Q}(x_\im^+,0) = -\partial_x \mathbf{F}^{(a)}(\boldsymbol{p}_i(x_\im)) - \mathbf{P}(\boldsymbol{p}_i(x_\im))\partial_{x}\mathbf{Q}(x_\im^+,0),\\
 \partial_{t}\mathbf{Q}(x_\ip^-,0) = -\partial_x \mathbf{F}^{(a)}(\boldsymbol{p}_i(x_\ip)) - \mathbf{P}(\boldsymbol{p}_i(x_\ip))\partial_{x}\mathbf{Q}(x_\ip^-,0).
\end{array}
\right\}
\end{equation}
The advection flux gradient, to second-order, can be approximated as follows
\begin{equation} \label{eq:flux_grad}
\partial_x \mathbf{F}^{(a)}(\boldsymbol{p}_i(x_\im))=\partial_x \mathbf{F}^{(a)}(\boldsymbol{p}_i(x_\ip))= \frac{\mathbf{F}^{(a)}(\mathbf{Q}_i^R)-\mathbf{F}^{(a)}(\mathbf{Q}_i^L) }{\Delta x}
\end{equation}
and the pressure non-conservative term as 
\begin{equation} \label{eq:flux_grad_1}
\left.
\begin{array}{ll}
\mathbf{P}(\boldsymbol{p}_i(x_\im))\partial_{x}\mathbf{Q}(x_\im^+,0) = \mathbf{P}(\mathbf{Q}_i^L)\mathbf{\Delta}_i,\\
\mathbf{P}(\boldsymbol{p}_i(x_\ip))\partial_{x}\mathbf{Q}(x_\ip^-,0) =  \mathbf{P}(\mathbf{Q}_i^R)\mathbf{\Delta}_i.
\end{array}
\right\}
\end{equation}
Finally, the evolved boundary values in cell $I_i$, at time $\tau = \frac{1}{2}\Delta t$, after using (\ref{eq:boundary_evolved}), (\ref{eq:flux_grad}) and (\ref{eq:flux_grad_1}), become
\begin{equation} \label{eq:extrapolated_values}
\left.
\begin{array}{ll}
\QevL = \mathbf{Q}_i^L - \um \Delta t \dfrac{\mathbf{F}^{(a)}(\mathbf{Q}_i^R) - \mathbf{F}^{(a)}(\mathbf{Q}_i^L) }{\Delta x} - \um \Delta t \mathbf{P}(\mathbf{Q}_i^L)\mathbf{\Delta}_i,\\
\QevR =  \mathbf{Q}_i^R - \um \Delta t \dfrac{\mathbf{F}^{(a)}(\mathbf{Q}_i^R) - \mathbf{F}^{(a)}(\mathbf{Q}_i^L) }{\Delta x}-\um \Delta t \mathbf{P}(\mathbf{Q}_i^R)\mathbf{\Delta}_i.
\end{array}
\right\}
\end{equation}
The time evolution is obtained in a splitting mode, making use of the advection fluxes $\mathbf{F}^{(a)}$ and the pressure coefficient matrix $\mathbf{P}$. We remark that we do not make use of the coefficient matrix of the full system $\mathbf{A}$ (\ref{eq:coeff_matrix}).

\subsection{The Fully Discrete Second Order Accurate One-Step Scheme}
\label{sec.rec.fullydiscrete}
Exact integration of the 
system (\ref{completesystem_cons_non_cons_no_source}) over a space-time control volume $V_i=[x_{i-\frac{1}{2}};x_{i+\frac{1}{2}}] \times [t^n;t^{n+1}]$ (see \cite{Castro2006} and \cite{Pares2006} for details) gives the following update formula:
\begin{equation} \label{eq:update_2nd_order}
\mathbf{Q}_i^{n+1} = \mathbf{Q}_i^{n} -  \frac{\Delta t}{\Delta x}\left[ \left( \mathbf{D}_{\ip}^- + \mathbf{D}_{\im}^+\right) + \left(\mathbf{F}^{(a)}_{\ip} -\mathbf{F}^{(a)}_{\im} \right)\right] -\Delta t \mathbf{H}_i     \;
\end{equation} 
where 
\begin{equation} \label{eq:Hi}
\mathbf{H}_i = \frac{1}{\Delta t} \displaystyle\int_0^{\Delta t} \displaystyle\int_0^1 \mathbf{P}(\mathbf{Q}(x_i,\um \Delta t)) \partial_x \mathbf{Q}(x_i,\um \Delta t) dx \, dt .
\end{equation}
The term $\mathbf{H}_i$ integrates the smooth part of the non-conservative product within each cell (excluding the jumps at the boundaries) and vanishes for a first order scheme where we have $\mathbf{H}_i=0$. This term will be defined shortly. 

\begin{figure}[tbp]
\centering
\includegraphics[width=0.4\columnwidth]{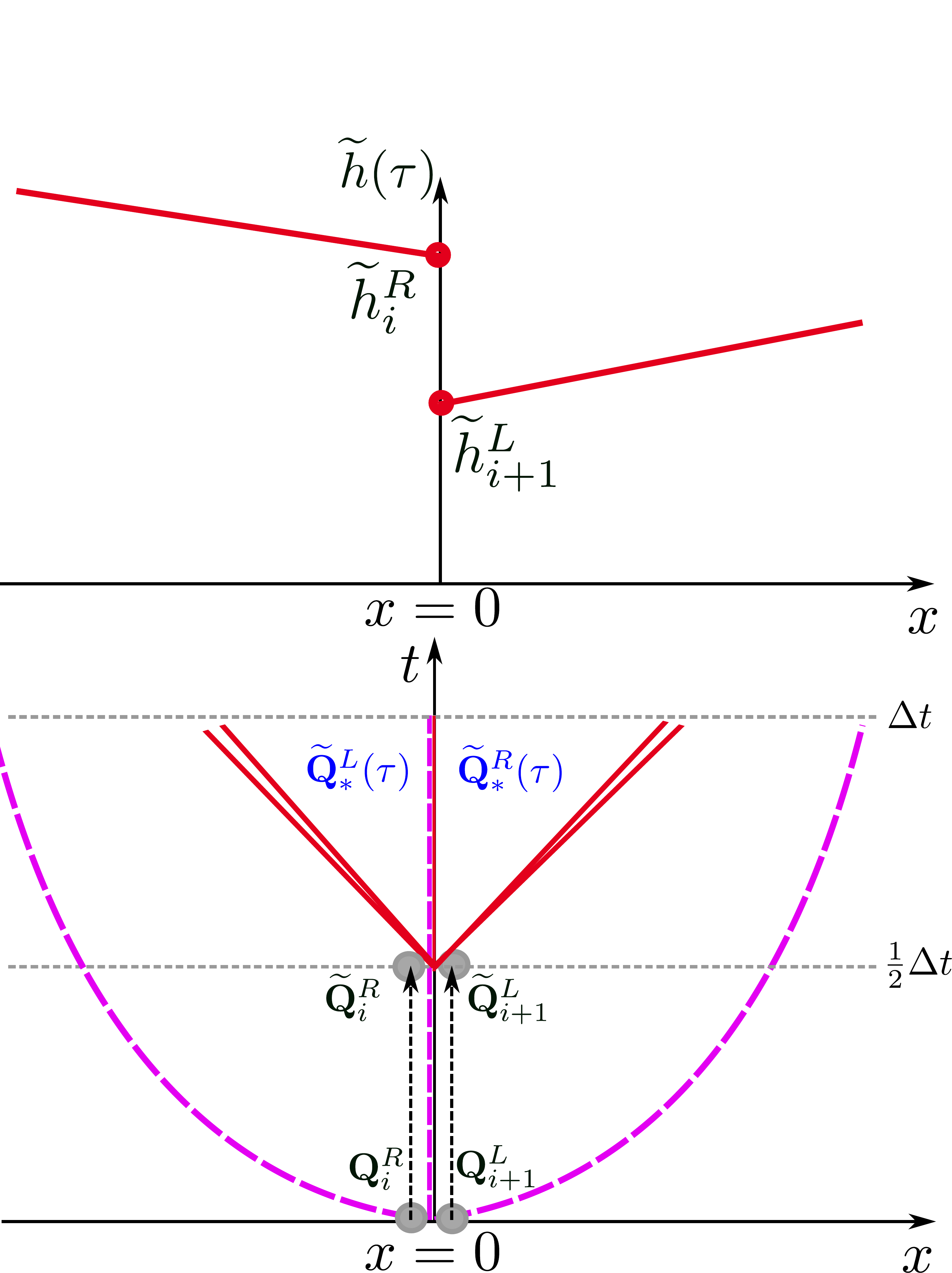}
\caption{ \textbf{Solution of GRP by HEOC method}.
\textbf{Top panel}: initial condition for a single
component ($\widetilde{h}$ in this example) of the vector of the evolved boundary values.
\textbf{Bottom panel}: the structure of the solution of the GRP is represented by the curved characteristics; the evolution step for half a time step is represented by the straight vertical arrows either side of the interface, and the structure of the solution of the conventional Riemann problem placed at the half time is represented by straight red characteristics. The sought solution at the fixed interface position $x_{\ip}$ , or $x=0$ in local coordinates, are $\widetilde{\mathbf{Q}}^{R}_*$ right at the interface and $\widetilde{\mathbf{Q}}^{L}_*$ on the left.
}
\label{fig:GRP}
\end{figure}

Now, the required solution of the generalized Riemann problem at the interface $x_\ip$ is given by the solution of the following conventional (piece-wise-constant data, homogeneous) Riemann problem (see Fig.~\ref{fig:GRP})
\begin{equation}                                                 \label{eq:GRP_prex}
       \left.
       \begin{array}{l}
       \partial_{t}\mathbf{Q} + \mathbf{P}(\mathbf{Q}) \partial_{x}\mathbf{Q} = 
       {\bf 0} \;,\hspace{2mm} x \in {\cal{R}}  \;, \hspace{2mm} t>0 \;, \\
       \\
       \mathbf{Q}(x,0) = \left\{ 
                      \begin{array}{lll}
                            {\widetilde{\mathbf{Q}}_i^R} & \mbox{ if } & x < 0  \;,\\
                            \\
                            {\widetilde{\mathbf{Q}}_{i+1}^L} & \mbox{ if } & x > 0  \;.
                      \end{array}
                      \right.
      \end{array}
      \right\}
\end{equation}
The solutions of the GRP (\ref{eq:GRP_prex})  are obtained solving the non-linear system (\ref{eq_non_lin_sys_h1}-\ref{eq_non_lin_sys_q}) as described in section \ref{sec:pressure_system} with initial data (${{\widetilde{\mathbf{Q}}_i^R},{\widetilde{\mathbf{Q}}_{i+1}^L}}$).
Let us denote the similarity solution of (\ref{eq:GRP_prex}) as $\widetilde{\mathbf{Q}}^{L}_*=[\tilde{h}_{*L},\tilde{q}_{*},\tilde{\eta}_{L}]^T$ and $\widetilde{\mathbf{Q}}^{R}_*=[\tilde{h}_{*R},\tilde{q}_{*},\tilde{\eta}_{R}]^T$, then 
the fluctuations $\mathbf{D}_{\ip}^-$ and $\mathbf{D}_{\ip}^+$ are calculated as 
\begin{equation} \label{eq:fluctuations_God_2nd_order}
\left.
\begin{array}{ll}
\mathbf{D}_{\ip}^- = \frac{1}{\Delta t} \displaystyle\int_{0}^{\Delta t} \displaystyle\int_{0}^{1} \mathbf{P}(\mathbf{\Psi}(s; {\widetilde{\mathbf{Q}}_i^R},\widetilde{\mathbf{Q}}_{*}^R)) \frac{ \partial}{\partial{s}}\mathbf{\Psi}(s; {\widetilde{\mathbf{Q}}_i^R},\widetilde{\mathbf{Q}}_{*}^R)) ds \, dt \\
\mathbf{D}_{\ip}^+ = \frac{1}{\Delta t} \displaystyle\int_{0}^{\Delta t} \displaystyle\int_{0}^{1} \mathbf{P}(\mathbf{\Psi}(s; \widetilde{\mathbf{Q}}^{*}_L,{\widetilde{\mathbf{Q}}_{i+1}^L})) \frac{ \partial}{\partial{s}}\mathbf{\Psi}(s; \widetilde{\mathbf{Q}}^{L}_*,{\widetilde{\mathbf{Q}}_{i+1}^L})) ds \, dt\;.
\end{array}
\right\}
\end{equation} 
Here, the paths considered are
\begin{equation}
\label{eq:seg_path_2nd}
\left.
\begin{array}{lll}
 \mathbf{\Psi}(s; {\widetilde{\mathbf{Q}}_i^R},\widetilde{\mathbf{Q}}_{R}^*) &= 
 \widetilde{\mathbf{Q}}_{i}^R 
 + s \left( \widetilde{\mathbf{Q}}_{R}^* - \widetilde{\mathbf{Q}}_i^R \right) \\ 
 \mathbf{\Psi}(s; {\widetilde{\mathbf{Q}}_{*}^L},\widetilde{\mathbf{Q}}_{i+1}^L) &= \widetilde{\mathbf{Q}}_{*}^L 
 + s \left(\widetilde{\mathbf{Q}}_{i+1}^L - \widetilde{\mathbf{Q}}_{*}^L  \right). 
\end{array}
\right\}  
\end{equation} 
 
Following the same approach as for the first order problem (see equation \ref{eq:advection_flux}), the numerical flux $\mathbf{F}^{(a)}_{\ip}$ are obtained as
\begin{equation}                         \label{eq:advection_flux_2nd}
     {\bf F}_{i+\frac{1}{2}}^{(a)} =  \left\{\begin{array}{c}
     \tilde{q}_{*} \begin{bmatrix}
     0\\
     \displaystyle{\tilde{u}_L} \\
     A_g \dfrac{\tilde{u}_L^{m-1}}{\tilde{h}_L} 
     \end{bmatrix}  \mbox{ if } \tilde{q}_{*} \ge 0 \;, \\
     \\
     \tilde{q}_{*} \begin{bmatrix}
     0\\
     \displaystyle{\tilde{u}_{i+1}} \\
     A_g \dfrac{\tilde{u}_{R}^{m-1}}{\tilde{h}_{R}} 
     \end{bmatrix}  \mbox{ if } \tilde{q}_{*} < 0 \;.
     \end{array} \right. 
\end{equation}
Finally, to compute $\mathbf{H}_i$ we first approximate the spatial derivative as
\begin{equation} \label{eq:spatial_deriv}
\partial_x \mathbf{Q}(x_i,\um \Delta t) =\frac{\QevR-\QevL}{\Delta x} .
\end{equation}
Substitution into equation (\ref{eq:Hi}) and integrating we obtain
\begin{equation} \label{eq:Hi_num}
\mathbf{H}_i =  \mathbf{P}(\mathbf{Q}(x_i,\um \Delta t)) \frac{\QevR-\QevL}{\Delta x} \,
\end{equation}
where
\begin{equation} \label{eq:Qi_evol}
\mathbf{Q}(x_i,\um \Delta t) = \mathbf{Q}_i^n  - \um \Delta t \mathbf{P}(\mathbf{Q}_i^n)\mathbf{\Delta}_i.
\end{equation}



In the following we briefly summarize the entire second-order one-step 
algorithm: 
\begin{enumerate}
 \item Perform the AENO reconstruction described in section \ref{sec.rec.nonlinear} in order to obtain the slope 
$\mathbf{\Delta}_i$ (\ref{eq:slope_AENO}) for each cell.
\item Extrapolate values at cell boundaries $x_\im$ and $x_\ip$ (\ref{eq:boundary_values}) and then evolve 
evolve these limiting values in time using the time Taylor series expansion (\ref{eq:boundary_evolved}).
 \item Solve the the GRP (\ref{eq:GRP_prex}) of the pressure system
through the solution of the non-linear system (\ref{eq_non_lin_sys_h1}-\ref{eq_non_lin_sys_q}) with  initial  data ($\widetilde{\mathbf{Q}}^{R}_i$, $\widetilde{\mathbf{Q}}^{L}_{i+1}$). This step gives the star region solutions $\widetilde{\mathbf{Q}}^{L}_*=[\tilde{h}_{*L},\tilde{q}_{*},\tilde{\eta}_{L}]^T$ and $\widetilde{\mathbf{Q}}_{*}^R=[\tilde{h}_{*R},\tilde{q}_{*},\tilde{\eta}_{R}]^T$.
 \item At each interface: use  ($\widetilde{\mathbf{Q}}^{R}_i$, $\widetilde{\mathbf{Q}}^{L}_{i+1}$) and ($\widetilde{\mathbf{Q}}_{*}^L$, $\widetilde{\mathbf{Q}}_{*}^R$) 
 to compute the fluctuations ( \ref{eq:fluctuations_God_2nd_order}) using the paths (\ref{eq:seg_path_2nd}).
\item  At each interface: use  ($\widetilde{\mathbf{Q}}^{R}_i$, $\widetilde{\mathbf{Q}}^{L}_{i+1}$) and ($\widetilde{\mathbf{Q}}_{*}^L$, $\widetilde{\mathbf{Q}}_{*}^R$) 
 to calculate the fluxes (\ref{eq:advection_flux_2nd}).
  \item At each cell center: use  ($\widetilde{\mathbf{Q}}^{R}_i$, $\widetilde{\mathbf{Q}}^{L}_{i+1}$) and time evolution of the center cell values (\ref{eq:Qi_evol}) to calculate $\mathbf{H}_i$ (\ref{eq:Hi_num}). 
 \item Finally use the fully discrete scheme (\ref{eq:update_2nd_order}) and perform the update of the cell averages. 
\end{enumerate}

\section{Numerical results}
\label{sec:num_res}
Here we assess the proposed splitting method on a carefully selected, suite of test problems. For all tests the numerical stability is imposed by the Courant-Friedrichs-Lewy condition 
and the integration time step is evaluated as
\begin{equation}
\Delta t = CFL \min_{1 \leq i \leq M} \frac{\Delta x}{\lambda_{Hi}} 
\end{equation}
where $M$ is the total number of cells and $\lambda_{Hi}=|q_i|/h_i + \sqrt{gh_i}$ is the maximum eigenvalue for the fixed bed case (Saint-Venant equations). To take into account the small differences between $\lambda_{Hi}$ and the maximum eigenvalue of the coupled SVE model (see \cite{Lyn:2002} for details) we set the $CFL$ number to 0.9 for all numerical runs.
As numerical reference schemes we use: (i) the Dumbser-Osher-Toro
(DOT) solver \cite{dumbser:2011}, which is an all-purpose universal Godunov upwind method, that can be applied to any hyperbolic system, as long as the full eigenstructure
is available; (ii) the PRICE-C scheme \cite{Canestrelli:2009}, which is a method of the centred type which requires a minimum knowledge about the eigenstructure, i.e. an estimate of the fastest eigenvalue to be used in the CFL condition.

\subsection{Verification of the C-property}
A desirable feature of numerical methods for shallow water systems with variable
bottom is the satisfaction of  the so-called C-property as introduced by \citet{Bermudez:1994}. Let us consider a quiescent flow ($q$=\SI{0}{m^2/s^{-1}}) over any submerged bed profile. Under these conditions the initial
water surface $H =h + \eta$ is constant and should remain constant in time. This is numerically achieved if, 
the solution does not change in time and thus $\mathbf{Q}_i^{n+1}=\mathbf{Q}_i^{n}$ in (\ref{eq:update}). Therefore we have to prove that  
\begin{equation} \label{eq_relat_C_property}
 \left(\mathbf{D}_{\im}^+  + \mathbf{D}_{\ip}^- \right) + \left(\mathbf{F}^{(a)}_{\ip} -\mathbf{F}^{(a)}_{\im} \right) =0   \;.
\end{equation}

First we consider the flux $\mathbf{F}^{(a)}_{\ip}$  computed as in (\ref{eq:advection_flux}). Under quiescent flow conditions  $q_L=q_R=0$ and consequently from ~\eqref{eq:K_delta_eta} we obtain $K= h_L^{3/2} + h_R^{3/2}$. Substitution of this result in (\ref{eq_non_lin_sys_h1}) gives 
\begin{equation} \label{eq:C_property_proof}
h_L^{3/2} -h_R^{3/2} - h_{*L}^{3/2} + h_{*R}^{3/2} = 0
\end{equation}
and thus we have from (\ref{eq_non_lin_sys_q}) $q_*=0$ . Inserting $q_*=0$  into (\ref{eq:advection_flux}) leads to $\mathbf{F}^{(a)}_{\ip}=0$. Analogous conclusions can be drawn  for $\mathbf{F}^{(a)}_{\im}$.    
    
Second, we focus on the fluctuations $\mathbf{D}_{\ip}^+$. The numerical evaluation of such fluctuation through (\ref{eq:fluctuations_God_num}) gives 
\begin{equation} 
\mathbf{D}_{\ip}^- = \mathbf{\widehat{P}}_\ip^+ \left[\mathbf{Q}^n_{i+1} -\mathbf{Q}_{*}^L  \right] =
\left[
\begin{array}{ccc}
0 & 1& 0 \\
g \bar{h} & 0 & g \bar{h} \\
0& 0 & 0 
\end{array}
\right]\; 
\left[ \begin{array}{c} 
 h_R-h_{*L} \\
  0\\
 \eta_R -\eta_L 
\end{array} \right]=
\left[ \begin{array}{c} 
 0 \\
  g  \bar{h} \left[ (h_R + \eta_R) - (h_{*L} + \eta_L)\right] \\
 0
\end{array} \right]
\end{equation}
with $\bar{h}=\int \limits_0^1 h(s) ds = \int \limits_0^1 ( h_L + s (h_{*R} - h_L)) ds$. 
If $q=0$, then (\ref{eq:Riemann_inv_left}) gives that $h_{*L}= h_L $ and thus the second element of  $\mathbf{D}_{\ip}^-$ becomes 
\begin{equation}
g  \bar{h} \left[ (h_R + \eta_R) - (h_{*L} + \eta_L)\right]=g  \bar{h} \left[ (h_{R} + \eta_R) - (h_{L} + \eta_L)\right] = g  \bar{h} (H_R-H_L) \;.
\end{equation}
Since $H$ is constant, also $H_L=H_R=H$ and therefore  the second element of $\mathbf{D}_{\ip}^-$ is zero. Thus 
we obtain that $\mathbf{D}_{\ip}^-=0$. 
The proof that  $\mathbf{D}_{\im}^+=0$  is found in
an entirely analogous way. 
Therefore all four terms in (\ref{eq_relat_C_property}) are identically zero and this demonstrates 
that our first-order splitting scheme is exactly well-balanced.
Finally, we remark that numerical tests conducted with our second order extension of the splitting method demonstrate that the scheme we propose is also well balanced (results not shown). 


\subsection{Numerical convergence study}
Here we verify the accuracy of our numerical scheme by studying empirical convergence rates  using the
method of manufactured solutions. 
For the assessment we compare the exact solutions against the numerical solutions employing AENO and ENO reconstruction.
We omit the presentation of ENO reconstruction as
the reader can consult Chapter 20 of \cite{Toro:2013}. 
We use the manufactured solutions presented in \cite{Canestrelli:2009} and proceed as follows. 
We consider the frictionless SVE equations
\begin{eqnarray}\label{system_convergence}
\left\{
\begin{array}{lc}
 \partial_t h  + \partial_x q &=0, \\ 
\partial_t q + \partial_x
\left( qu+\frac{1}{2}gh^2\right)+ gh \partial_x \eta&=0, \\
\partial_t \eta + \partial_x q_b &= 0 
\end{array} 
\right.
\end{eqnarray}
and prescribe 
three smooth functions for $h(x,t)$, $q(x,t)$ and $\eta(x,t)$ which 
satisfy exactly (\ref{system_convergence}). These functions are 
\begin{eqnarray}
\left\{
\begin{array}{ll}
h(x,t) & = h_0 +c_0 \sin (k x -\omega t), \\ 
q(x,t) &= \frac{\omega}{k} h_0 + c_0 \frac{\omega}{k} \sin (k
-\omega t), \\
\eta(x,t) & = -h(x,t),\\
q_b(x,t) &= -q(x,t)  
\end{array}
\right. 
\end{eqnarray}
with $k = 2 \pi/L_w$ and $\omega = 2 \pi /T_p$.

\begin{table} [t]
\caption{\textbf{Convergence-rate study for the sediment transport problem}. Splitting-ADER scheme
for the 2\textsuperscript{nd} order of accuracy with the \textbf{AENO} reconstruction. Rates are
calculated for the discharge per unit width $q$ and bed level $\eta$. Sediment transport is quantified by using the sediment transport formula \eqref{eq:Grass} with $A_g=0.01$ and $m=1.5$. Computational parameters are: $T_{final}$=\SI{10}{s}, domain length $L$=\SI{500}{m},
$CFL$=0.9, $c_0$=\SI{0.01}{m}, $h_0$=\SI{5}{m}, $T_p$=\SI{10}{s}, $L_w$=\SI{250}{m}. AENO reconstruction is performed with $TOL=10^{-4}$ and $\epsilon$=1.}
\label{tab:conv_AENO}
\centering
\small
\renewcommand{\arraystretch}{1.0}
\begin{tabular}{c | c | c | c | c | c | c | c | c | c}
\hline\hline
\multicolumn{1}{c|}{} & \multicolumn{4}{c|}{variable $q$} & \multicolumn{4}{c|}{variable ${\eta}$} & \multicolumn{1}{c}{} \\
\hline
M & $L_1$ & $\mathcal{O}(L_1)$ & $L_{\infty}$ & $\mathcal{O}(L_\infty)$ &
    $L_1$ & $\mathcal{O}(L_1)$ & $L_{\infty}$ & $\mathcal{O}(L_\infty)$ & CPU [$s$] \\
\hline
    20&  2.15E-04&    -   &  4.01E-04&    -   &  2.45E-05&    -   &  4.82E-05&  -   & 0.05\\
    40&  5.31E-05&    2.02&  7.21E-05&    2.48&  6.57E-06&    1.90&  1.85E-05&  1.38& 0.07\\
    80&  1.24E-05&    2.10&  1.56E-05&    2.21&  1.65E-06&    1.99&  5.55E-06&  1.73& 0.14\\
   160&  3.03E-06&    2.03&  3.67E-06&    2.09&  3.99E-06&    2.05&  1.89E-06&  1.56& 0.45\\
   320&  7.48E-07&    2.01&  8.91E-07&    2.04&  9.93E-08&    2.01&  6.51E-07&  1.54& 2.04\\
   640&  1.86E-07&    2.01&  2.20E-07&    2.02&  2.47E-08&    2.00&  2.31E-07&  1.49& 7.20\\
   1280& 4.65E-08&    2.00&  5.46E-08&    2.01&  6.23E-09&    1.99&  8.10E-08&  1.51& 29.37\\
\hline\hline
\end{tabular}
\normalsize
\end{table}

\begin{table} [tbp]
\caption{\textbf{Convergence-rate study for the sediment transport problem}. Splitting-ADER scheme
for the 2\textsuperscript{nd} order of accuracy with the \textbf{ENO} reconstruction. Rates are
calculated for the discharge per unit width $q$ and bed level $\eta$. Sediment transport is quantified by using the sediment transport formula \eqref{eq:Grass} with $A_g=0.01$ and $m=1.5$.
Computational parameters are: $T_{final}$=\SI{10}{s}, domain length $L$=\SI{500}{m},
$CFL$=0.9, $c_0$=\SI{0.01}{m}, $h_0$=\SI{5}{m}, $T_p$=\SI{10}{s}, $L_w$=\SI{250}{m}.}
\label{tab:conv_ENO}
\centering
\small
\renewcommand{\arraystretch}{1.0}
\begin{tabular}{c | c | c | c | c | c | c | c | c | c}
\hline\hline
\multicolumn{1}{c|}{} & \multicolumn{4}{c|}{variable $q$} & \multicolumn{4}{c|}{variable $\eta$} & \multicolumn{1}{c}{} \\
\hline
M & $L_1$ & $\mathcal{O}(L_1)$ & $L_{\infty}$ & $\mathcal{O}(L_\infty)$ &
    $L_1$ & $\mathcal{O}(L_1)$ & $L_{\infty}$ & $\mathcal{O}(L_\infty)$ & CPU [$s$] \\
\hline
    20&  2.89E-04&    -   &  6.26E-04&    -   &  7.62E-05&    -   &  1.32E-04&  -   & 0.04 \\
    40&  1.05E-04&    1.46&  2.66E-04&    1.24&  4.91E-05&    0.64&  8.00E-05&  0.73& 0.07 \\
    80&  3.44E-05&    1.61&  1.12E-04&    1.25&  1.70E-05&    1.52&  5.40E-05&  0.57& 0.13 \\
   160&  9.89E-06&    1.80&  4.49E-05&    1.31&  6.43E-06&    1.41&  3.56E-05&  0.60& 0.47 \\
   320&  2.70E-06&    1.87&  1.90E-05&    1.24&  2.66E-06&    1.27&  2.46E-05&  0.54& 1.81 \\
   640&  7.92E-07&    1.77&  8.04E-06&    1.24&  1.18E-06&    1.18&  1.82E-05&  0.44& 7.09 \\
   1280&  2.63E-07&    1.59&  3.42E-06&    1.23& 6.34E-07&    0.89&  1.46E-05&  0.30& 29.47 \\
\hline\hline
\end{tabular}
\normalsize
\end{table}
  
Results are presented in terms of standard
norms $L_1$, $L_{\infty}$  and relative convergence rates for variables $q$ and $\eta$ and given in Tables \ref{tab:conv_AENO} and \ref{tab:conv_ENO}. As expected from the settings of this test, the results for the variable $h$ are similar to those of $\eta$ and are not reported here. 
The AENO method reached the expected rate in all norms, being
suboptimal in $L_{\infty}$ norm for $\eta$, similarly to what obtained in \cite{Toro:2020} in their application on blood flows. Although converging to the correct solution, the ENO method did not reach  the expected rate in all norms. From the comparison of the CPU time, the AENO and ENO reconstructions are comparable.

\subsection{Riemann problem test with movable and fixed bed}
\begin{figure}[tbp]
\centering
\includegraphics[width=0.9\columnwidth]{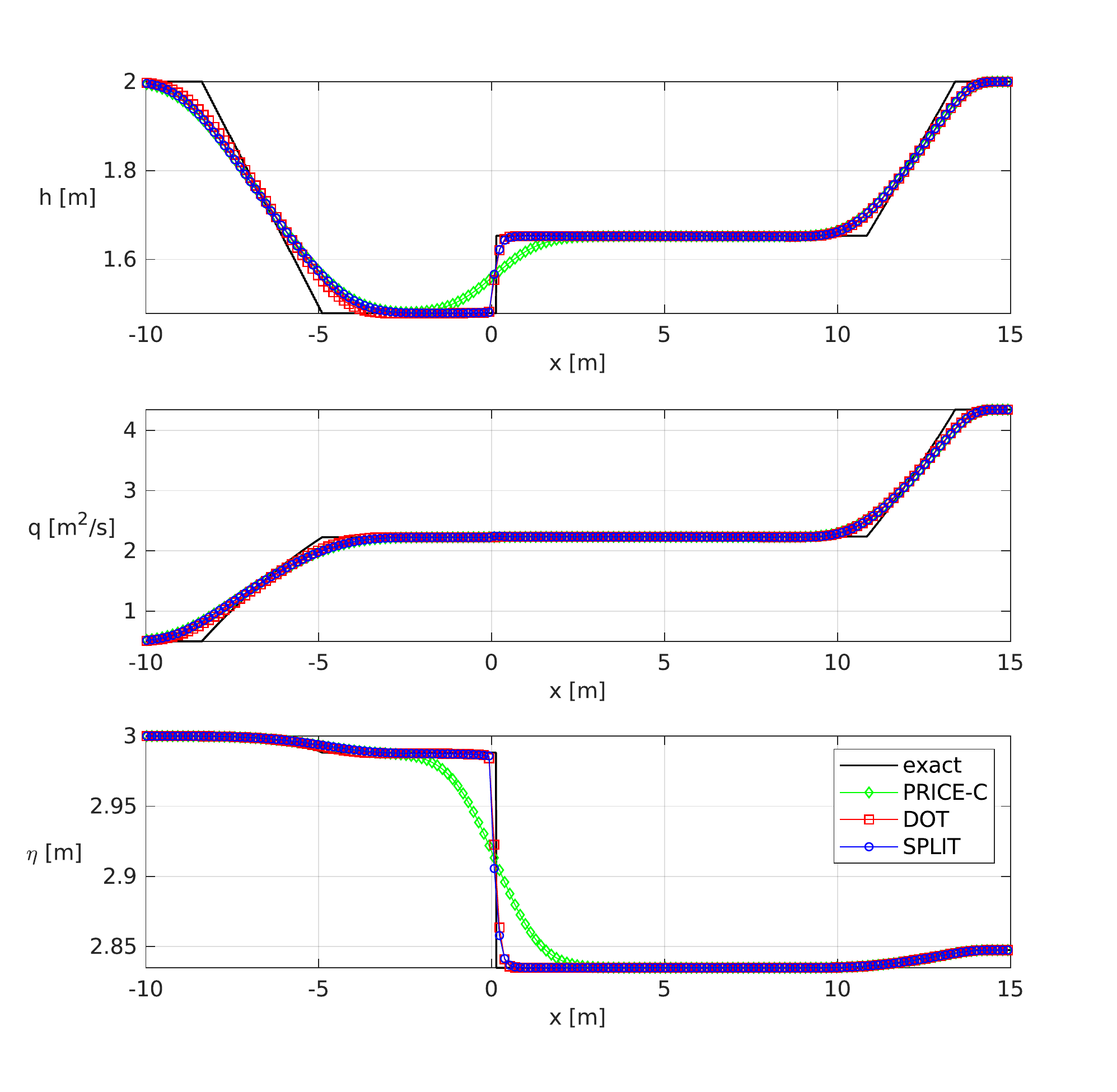}
\caption{\textbf{Results for a Riemann problem with with movable bed (1\textsuperscript{st} order of accuracy)}. 
Sediment transport is quantified by using the sediment transport formula \eqref{eq:Grass} with $A_g$=0.01 and $m$=3.0. 
Computational parameters are: $M$=200, domain length $L$=\SI{30}{m}, $T_{final}$=\SI{2}{s} and $CFL$=0.9.
}
\label{fig:T4_1ord}
\end{figure}

Here we assess the methods as applied to two  Riemann problem tests with exact
solution, one with movable bed \cite{murillo:2010} and one with fixed bed. 
The initial discontinuity is at $x$ = 0 and initial data  to the left and right are are given in Table ~\ref{tab:initial_RP}.
For the test on movable bed, the splitting numerical solutions are compared with respect to the PRICE-C and DOT at first order of accuracy and to DOT 
with AENO reconstruction in conjunction with ADER at second order. 
Fig.~\ref{fig:T4_1ord} shows results for the first order and demonstrate as in spite of its simplicity, results with our splitting are comparable with that of the
more sophisticated upwind DOT scheme. As expected the PRICE-C method considerably diffuses the central shock wave.   
\begin{table}
\caption{Initial conditions for Riemann problem tests.}
   \begin{center}
   \begin{tabular}{c|c|c|c|c|c|c}
   \hline\hline
       test  & $h_L$ [m]   &   $q_L$   [\si{m^2s^{-1}}] & $\eta_L$ [m] &  $h_R$ [m]   &   $q_R$   [\si{m^2s^{-1}}] & $\eta_R$ [m] \\ [3pt]    
       \hline    
        movable bed & 2.0 &  0.5 & 3.0  &  2.0  & 4.34297 & 2.84751 \\
        fixed bed & 1.0 &  0.0 & 0.0  &  0.1  & 0.0 & 0.0 \\
        \hline\hline
   \end{tabular}
   \label{tab:initial_RP}
   \end{center}
 \end{table}

\begin{figure}[tbp]
\centering
\includegraphics[width=0.9\columnwidth]{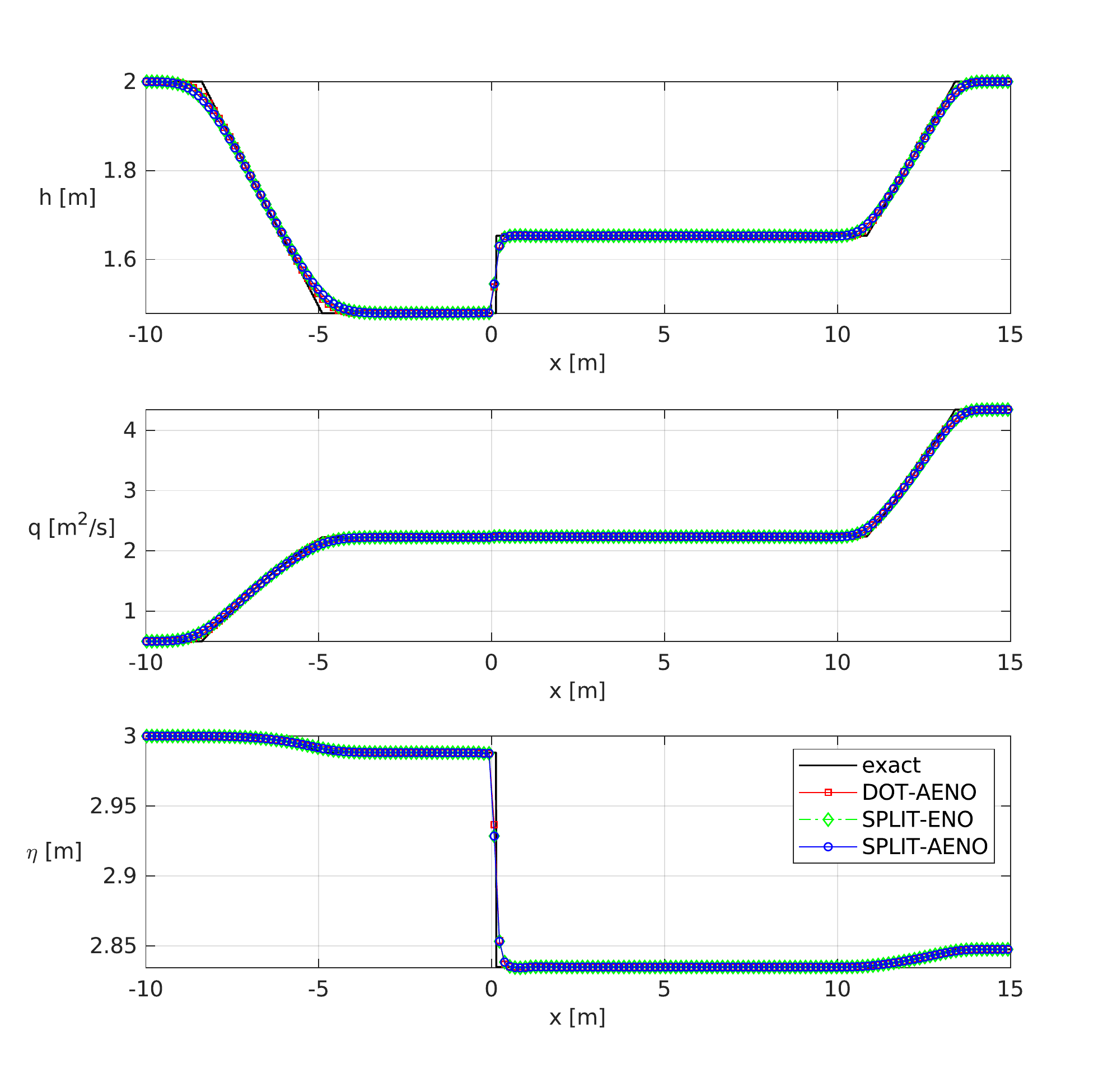}
\caption{\textbf{Results for a Riemann problem with movable bed (2\textsuperscript{nd} order of accuracy)}. 
Numerical solutions with the splitting and the DOT method are compared with the exact solution. 
Sediment transport is quantified by using the sediment transport formula \eqref{eq:Grass} with $A_g$=0.01 and $m$=3.0. 
Computational parameters are: $M$=200, domain length $L$=\SI{30}{m}, $T_{final}$=\SI{2}{s} and $CFL$=0.9. AENO reconstruction is performed using $TOL$= 10$^{-4}$ and $\epsilon$=0.5. 
}
\label{fig:T4_2ord}
\end{figure}
Results at second order of accuracy in Fig.~\ref{fig:T4_2ord} show as all three schemes give good results for this test problem, in
which the central shock wave moves very slowly.

\begin{figure}[tbp]
\centering
\includegraphics[width=0.9\columnwidth]{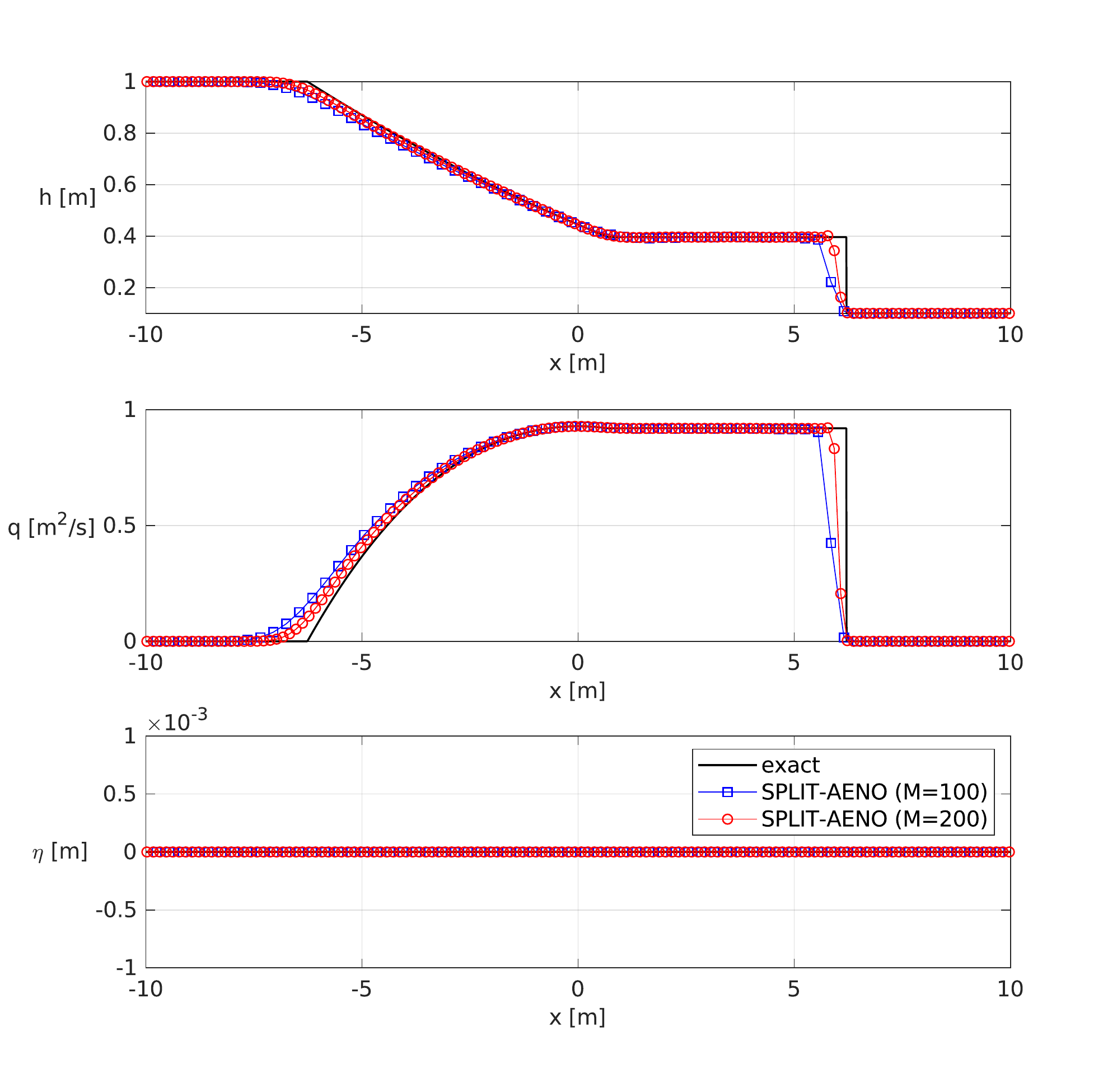}
\caption{\textbf{Results for a Riemann problem with fixed bed (2\textsuperscript{nd} order of accuracy)}.
Sediment transport is inhibited setting $A_g$=0 in the sediment transport formula \eqref{eq:Grass}.
Computational parameters are: $M$=100, domain length $L$=\SI{30}{m}, $T_{final}$=\SI{2}{s} and $CFL$=0.9.
AENO reconstruction is performed using $TOL$= 10$^{-4}$ and $\epsilon$=0.5.
}
\label{fig:T5_2ord}
\end{figure}

The results for the fixed bed case are displayed in Fig.~\ref{fig:T5_2ord}. They demonstrate an important feature of the proposed splitting method. That is, our method converges to the solution of the hydrodynamic  Saint-Venant equations when the flow does not transport sediments. This is a desirable feature that makes the proposed method particularly suitable for engineering applications.

The solution is composed by two external rarefaction waves and a central slowly moving shock. Our splitting method structure well describe the shock, both in terms of strength and position in all variables  providing very similar results to the DOT scheme also for the two rarefaction waves.

\subsection{Evolution of a sediment hump}
\begin{figure}[tbp]   
\centering
\includegraphics[width=0.9\columnwidth]{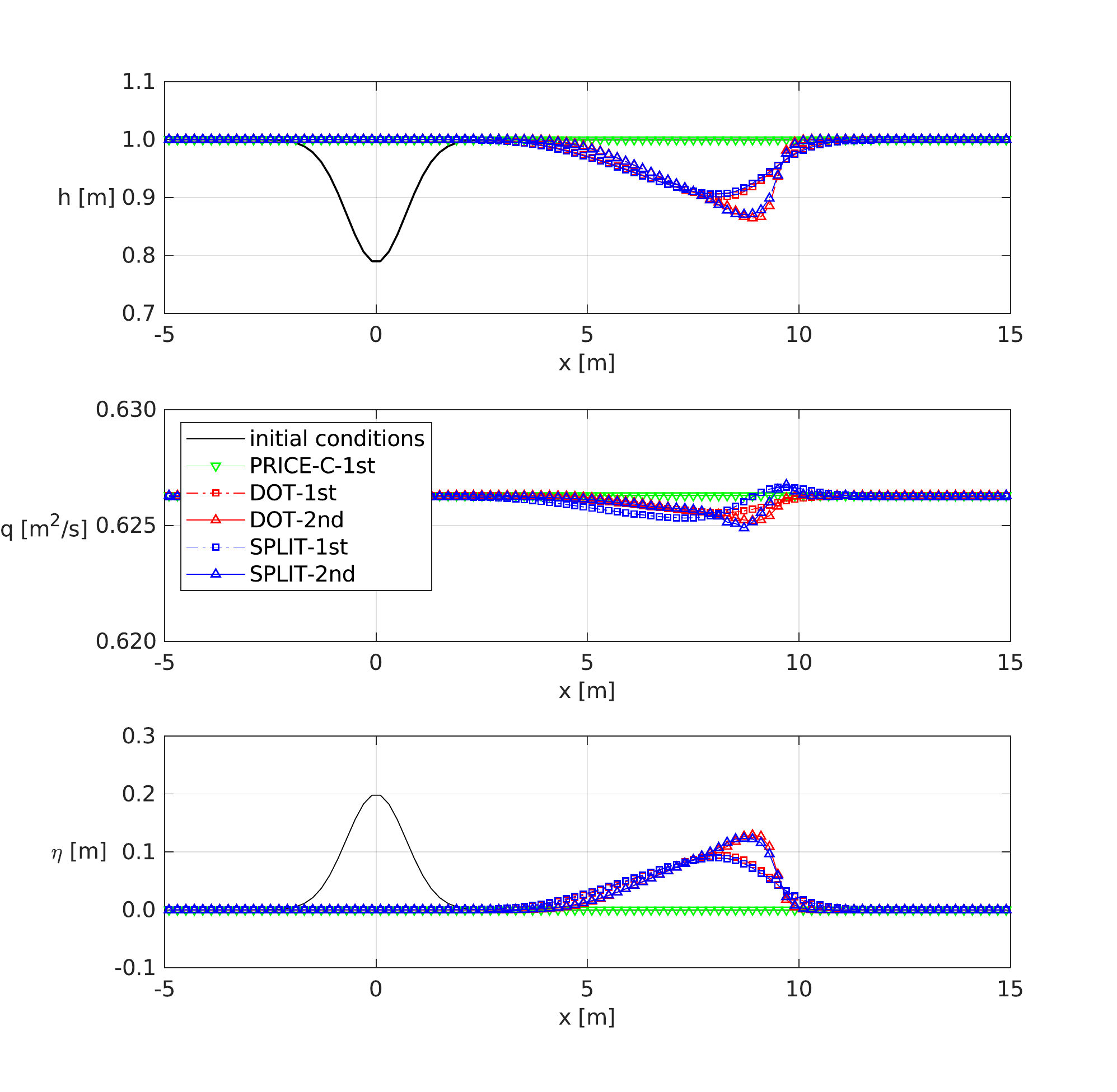}
\caption{\textbf{Long term evolution of a sediment hump under subcritical conditions ($Fr$=0.2)}. Numerical solutions of 1\textsuperscript{st} order (PRICE-C, DOT and splitting) and 2\textsuperscript{nd} order (splitting and DOT ADER with AENO reconstruction) are compared with the exact solution.
Sediment transport is quantified by using the sediment transport formula \eqref{eq:Grass_ucr} with $A_g$=0.01, $m$=1.5 and $\psi_u=10^{-3}$.
Computational parameters are: mesh $M$=100 cells, domain length $L$=\SI{20}{m} ($-10 \:\text{m} \leq x \leq 10 \:\text{m}$), $T_{final}$=\SI{1000}{s} and $CFL$=0.9. AENO reconstruction is performed using $TOL$=$10^{-4}$ and $\epsilon$=0.5. 
}
\label{fig:large_long_hump}
\end{figure}

This test case simulates the long term evolution of an erodible bed hump immersed into a quasi-steady, frictionless flow. 
The initial bed shape is described as   
\begin{equation} \label{eq:hump}
\eta(x,0)= \eta_{\max} \, e^{-x^2}  	,
\end{equation}
where $\eta_{\max}$  is the initial maximum hump amplitude set to $\eta_{\max}=0.2$ m.
For this test case the sediment transport formula is 
\begin{equation}  \label{eq:Grass_ucr}
	q_b = A_g \, [\max(u-u_{cr}) , 0]^m ,\qquad\text{with} \qquad u_{cr} = |u|- \left(\frac{\psi_u h_0}{m A_g}\right)^{\frac{1}{m-1}}
\end{equation}
where $u_{cr}$ is the fluid velocity critical value that must be exceeded for bedload transport to occur and $\psi_u$ is a small constant parameter. 
Flow discharge is kept constant at the inflow boundary and set equal to $q(x=0,t)$=\SI{0.6263}{m^2s^{-1}} while 
a constant water depth $h_0$=\SI{1}{m} is set at the downstream end of the domain. Transmissive downstream boundary conditions are set for the bed. Initial conditions correspond to the backwater profile obtained with these two boundary conditions.
The comparison is made for both 1\textsuperscript{st} and 2\textsuperscript{nd} order solutions. Results displayed in Fig.~\ref{fig:large_long_hump} demonstrate that the proposed scheme describes the hump evolution with good accuracy, very similar  to the one of the more sophisticated DOT method. This latter feature is particularly important when long term bed evolution must be studied.
Results also show that the PRICE-C scheme is extremely diffusive and the final time the hump is completed smeared.

\subsection{Short term propagation of a small sediment hump}
With this test we aim to reproduce bed movement under different flow conditions each characterized by a different Froude number. We consider a one-dimensional flat channel with a small hump on the bed described as  (\ref{eq:hump}).
For this test case we use the sediment transport formula (\ref{eq:Grass_ucr}).
We then consider  two different flow conditions, namely near-critical with $Fr=0.99$ and supercritical with $Fr=1.2$ and choose as reference water depth $h$=\SI{1}{m}. Given $Fr$ and $h$ we can calculate the value of the discharge per unit width $q$ which is kept always constant at the inlet. Thus, the initial condition $h(x,0)$ and $q(x,0)$ are obtained running the code under fixed bed conditions.   
We consider a very small hump and set $\eta_{\max}=10^{-5}$\si{m}. We then find the solution  by application of a linearized analytical solver, which is suitable for studying the
propagation of small-amplitude waves (see details in  \cite{Lyn:2002} and \cite{Canestrelli:2009}).  
\begin{figure}[tbp]
\centering
\includegraphics[width=0.9\columnwidth]{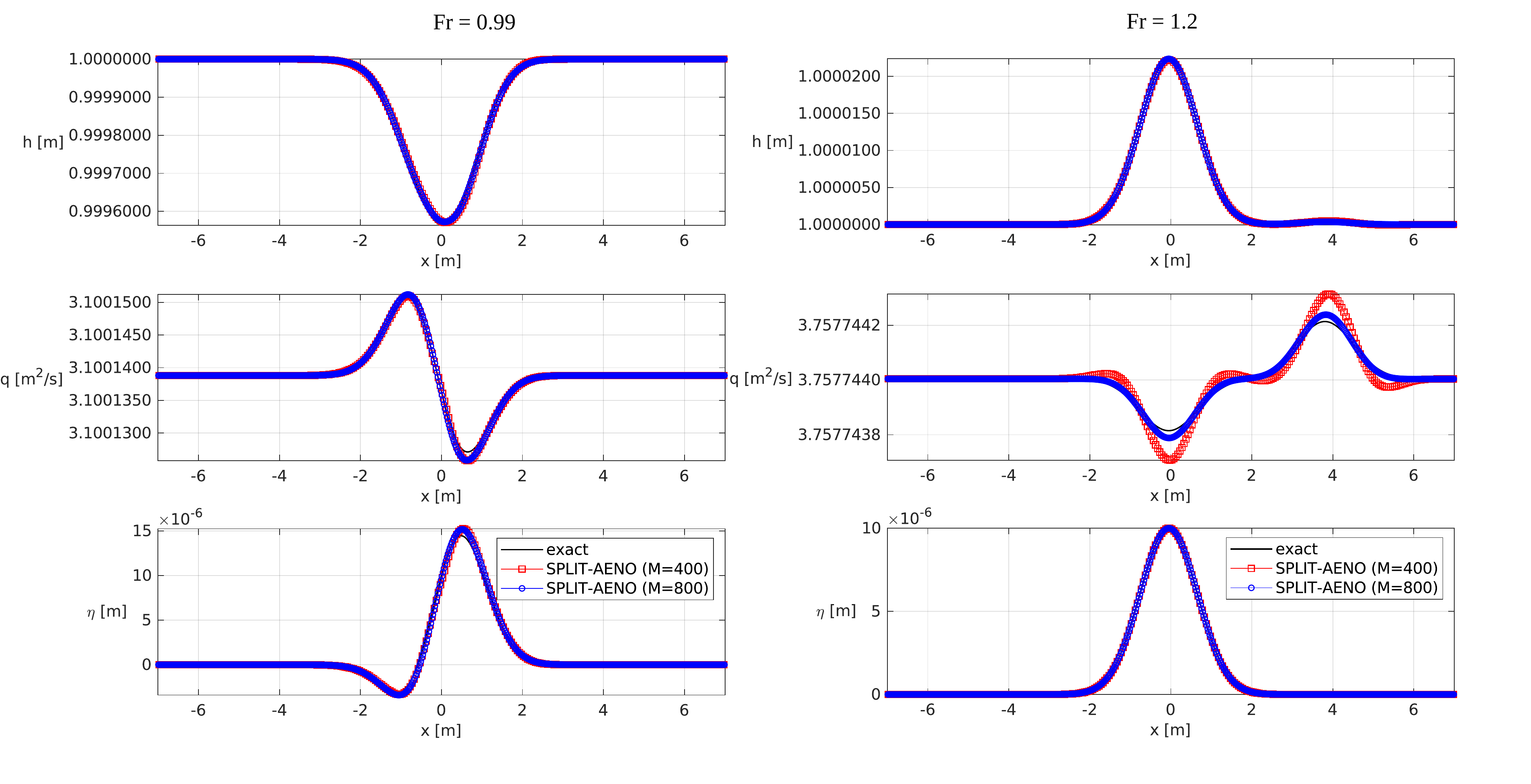}
\caption{\textbf{Short term propagation of a small sediment hump}.  \textbf{(Left panels}: near critical conditions ($Fr$=0.99)). \textbf{(Right panels}: supercritical conditions ($Fr$=1.2)).
Numerical solutions of 2\textsuperscript{nd} order splitting ADER with AENO reconstruction are compared with the linearized solution.
Sediment transport is quantified by using the sediment transport formula \eqref{eq:Grass_ucr} with $A_g$=0.01, $m$=1.5 and $\psi_u=10^{-2}$.
Computational parameters are: mesh $M$=[400, 800] cells, domain length $L$=\SI{20}{m}, $T_{final}$=\SI{6}{s} and $CFL$=0.9. AENO reconstruction is performed using $TOL$=$10^{-4}$ and $\epsilon$=0.5. 
}
\label{fig:small_hump}
\end{figure}
Numerical results are compared with the linearized solutions in Fig. (\ref{fig:small_hump}). In all cases they are in good agreement and converge to the linearized solution. In the supercritical case the numerical solution correctly predicts upstream propagation of the small hump.  In the transcritical case the two bed waves, one erosional  propagating upstream and one 
depositional propagating downstream are also correctly described.

\section{Concluding remarks}
\label{sec:conclusions}
In this paper we have proposed a splitting scheme for the SVE model which is valid in general for any hyperbolic  non-conservative system of PDEs. 
Then, for the SVE system we have studied the associated two systems of differential equations. After a careful study of the resulting two systems 
of PDEs we proposed a methodology for their numerical solution in the framework of Godunov methods. We then extend the the method up to the second order of accuracy. Finally we assess the robustness of
our splitting method considering different test cases. Results demonstrate that solutions converge correctly to second order
of accuracy in space and time, satisfy the well-balanced property  and are accurate when compared
with existing techniques. 
Our splitting method constitutes a building block for the construction of high-order  numerical methods and can easily include source terms at high order with the ADER approach. Furthermore, it gives the possibility of future extension to multiple space dimensions. 
The present approach also offers a simple way to incorporate  sediment transport formulas for the quantification of sediment fluxes.
These features are very attractive and makes the splitting scheme 
a viable alternative to existing approaches to be used for the
solution of river and near-shore engineering morphodynamic problems.

\section*{Acknowledgements}
{AS acknowledges financial support from the Italian Ministry of Education, University and Research (MIUR) via the Departments of Excellence initiative 2018–2022 attributed to DICAM of the University of Trento (grant L. 232/2016).
}


\end{document}